\documentclass[12pt,a4paper]{article}

\usepackage{amsmath}
\usepackage{amssymb}
\usepackage{amsthm}
\usepackage{t1enc}
\usepackage[latin1]{inputenc}
\usepackage[english]{babel}
\usepackage{amsfonts}
\usepackage{latexsym}
\usepackage{bm}
\usepackage{setspace}
\usepackage[toc,page]{appendix}
\usepackage{graphicx}
\usepackage{enumerate}
\usepackage[numbers,sort&compress]{natbib}
\usepackage{tikz}
\usepackage{csquotes}
\usepackage{doi}

\usepackage{subfigure}
\usepackage{floatrow}
\usepackage{mathtools}
\usepackage[font=footnotesize,labelfont=bf]{caption}

\DeclarePairedDelimiter{\floor}{\lfloor}{\rfloor}

\newtheorem{theorem}{Theorem}
\newtheorem{lemma}{Lemma}

\newtheorem{observation}{Observation}
\newtheorem{construction}{Construction}

\usepackage{color}   
\usepackage{hyperref}
\hypersetup{
    colorlinks=true, 
    linktoc=all,     
    linkcolor=blue,  
    hidelinks,
}

\newcommand{\cI}{{\cal I}}

\let\oldenumerate\enumerate
\renewcommand{\enumerate}{
  \oldenumerate
  \setlength{\itemsep}{0pt}
  \setlength{\parskip}{0pt}
  \setlength{\parsep}{0pt}
}

\renewenvironment{proof}{{\noindent\bfseries Proof. }}

\allowdisplaybreaks

\begin{document}


\title{Solution to a $3$-path isolation problem for subcubic graphs}
\author{Karl Bartolo\footnote{Email address: karl.bartolo.16@um.edu.mt} \quad \quad Peter Borg\footnote{Email address: peter.borg@um.edu.mt} \quad \quad Dayle Scicluna\footnote{Email address: dayle.scicluna.09@um.edu.mt} \\[5mm]
{\normalsize Department of Mathematics} \\
{\normalsize Faculty of Science} \\
{\normalsize University of Malta}\\
{\normalsize Malta}\\
}
\date{}

\maketitle

\begin{abstract}
The $3$-path isolation number of a connected $n$-vertex graph $G$, denoted by $\iota(G,P_3)$, is the size of a smallest subset $D$ of the vertex set of $G$ such that the closed neighbourhood $N[D]$ of $D$ in $G$ intersects the vertex sets of the $3$-vertex paths of $G$, meaning that no two edges of $G-N[D]$ intersect. If $G$ is not a $3$-path or a $3$-cycle or a $6$-cycle, then $\iota(G,P_3) \leq 2n/7$. This was proved by Zhang and Wu, and independently by Borg in a slightly extended form. The bound is attained by infinitely many connected graphs having induced $6$-cycles. Huang, Zhang and Jin showed that if $G$ has no $6$-cycles, or $G$ has no induced $5$-cycles and no induced $6$-cycles, then $\iota(G,P_3) \leq n/4$ unless $G$ is a $3$-path or a $3$-cycle or a $7$-cycle or an $11$-cycle. They asked if the bound still holds asymptotically for connected graphs having no induced $6$-cycles. Thus, the problem essentially is whether induced $6$-cycles solely account for the difference between the two bounds. In this paper, we solve this problem for subcubic graphs, which need to be treated differently from other graphs. We show that if $G$ is subcubic and has no induced $6$-cycles, then $\iota(G,P_3) \leq n/4$ unless $G$ is a copy of one of $12$ particular graphs whose orders are $3$, $7$, $11$ and $15$. The bound is sharp. 
\end{abstract}

\section{Introduction}

Unless explicitly stated, we use small letters to denote non-negative integers or elements of sets, and capital letters to denote sets or graphs. The set of positive integers is denoted by $\mathbb{N}$. For $n \geq 1$, $[n]$ denotes the set $\{i \in \mathbb{N} \colon i \leq n\}$. We take $[0]$ to be the empty set $\emptyset$. Arbitrary sets are taken to be finite. For a set $X$, $\binom{X}{2}$ denotes the set of $2$-element subsets of $X$ (that is, $\binom{X}{2} = \{ \{x,y \} \colon x,y \in X, x \neq y \}$). We may represent a $2$-element set $\{x,y\}$ by $xy$. For standard terminology in graph theory, the reader is referred to \cite{West}. Most of the terminology used in this paper is defined in \cite{Borg1}.

Every graph $G$ is taken to be \emph{simple}, that is, $G$ is a pair $(V(G),E(G))$ such that $V(G)$ (the vertex set of $G$) and $E(G)$ (the edge set of $G$) are sets satisfying $E(G) \subseteq \binom{V(G)}{2}$. If $|V(G)| = n$, then we call $G$ an $n$-\emph{vertex graph} or a graph of \emph{order} $n$. We call $G$ an \emph{$m$-edge graph} if $|E(G)| = m$. For each vertex $v$ of $G$, $N_G(v)$ denotes the set $\{w \in V(G) \colon vw \in E(G)\}$ of neighbours of $v$ in $G$, $N_G[v]$ denotes the closed neighbourhood $N_G(v) \cup \{ v \}$ of $v$, and $d_G(v)$ denotes the degree $|N_G (v)|$ of $v$. If $d_G(v) = 1$, then $v$ is called a \emph{leaf of $G$}. The minimum degree $\min\{d_G(v) \colon v \in V(G)\}$ and the maximum degree $\max\{d_G(v) \colon v \in V(G)\}$ of $G$ are denoted by $\delta(G)$ and $\Delta(G)$, respectively. If $\Delta(G) \leq 3$, then $G$ is said to be \emph{subcubic}. For a subset $X$ of $V(G)$, $N_G[X]$ denotes the closed neighbourhood $\bigcup_{v \in X} N_G[v]$ of $X$, $G[X]$ denotes the subgraph $(X,E(G) \cap \binom{X}{2})$ of $G$ induced by $X$, and $G - X$ denotes the graph $G[V(G) \setminus X]$ obtained by deleting the vertices in $X$ from $G$. We may abbreviate $G - \{x\}$ to $G-x$. For any $v,w \in V(G)$, the \emph{distance} between $v$ and $w$, denoted by $d_G(v, w)$, is the number of edges of a shortest path of $G$ from $v$ to $w$. Where no confusion arises, the subscript $G$ may be omitted; for example, $N_G(v)$ may be abbreviated to $N(v)$.

Consider two graphs $G$ and $H$. If $G$ is a copy of $H$, then we write $G \simeq H$ and we say that $G$ is an \emph{$H$-copy}. We say that \emph{$G$ contains $H$} if $H$ is a subgraph of $G$. We say that $H$ is an \emph{induced} subgraph of $G$ if $H = G[X]$ for some $X \subseteq V(G)$.

For $n \geq 1$, $K_n$ and $P_n$ denote the graphs $([n], \binom{[n]}{2})$ and $([n], \{\{i,i+1\} \colon i \in [n-1]\})$, respectively. For $n \geq 3$, $C_n$ denotes the graph $([n], E(P_n)\cup\{\{1,n\}\})$. A $K_n$-copy is called an \emph{$n$-clique} or a \emph{complete graph}, a $P_n$-copy is called an \emph{$n$-path} or simply a path, and a $C_n$-copy is called an \emph{$n$-cycle} or simply a cycle. Note that a $3$-clique is a $3$-cycle, and it is also called a \emph{triangle}. 

If $\mathcal{F}$ is a set of graphs and $F$ is a copy of a graph in $\mathcal{F}$, then we call $F$ an \emph{$\mathcal{F}$-graph}. If $G$ is a graph and $D \subseteq V(G)$ such that $N[D]$ intersects the vertex sets of the $\mathcal{F}$-graphs contained by $G$, then $D$ is called an \emph{$\mathcal{F}$-isolating set of $G$}. Thus, a subset $D$ of $V(G)$ is an $\mathcal{F}$-isolating set of $G$ if and only if $G - N[D]$ contains no $\mathcal{F}$-graph. The size of a smallest $\mathcal{F}$-isolating set of $G$ is denoted by $\iota(G, \mathcal{F})$ and called the \emph{$\mathcal{F}$-isolation number of $G$}. If $\mathcal{F} = \{F\}$, then we may replace $\mathcal{F}$ in the aforementioned terms and notation by $F$.

The study of isolating sets, initiated in the paper \cite{CaHa17} of Caro and Hansberg, is a natural generalization of the popular study of dominating sets \cite{C,CH,HHS,HHS2,HL,HL2}. Indeed, $D$ is a \emph{dominating set of $G$} (that is, $N[D] = V(G)$) if and only if $D$ is a $K_1$-isolating set of $G$, so the \emph{domination number of $G$}, denoted by $\gamma(G)$, is $\iota(G,K_1)$. Assume that $G$ is a connected $n$-vertex graph. One of the earliest domination results, due to Ore \cite{Ore}, is that $\gamma(G) \leq n/2$ if $G \not\simeq K_1$ (see \cite{HHS}). While the deletion of the closed neighbourhood of a dominating set gives the graph with no vertices, the deletion of the closed neighbourhood of a $K_2$-isolating set gives a graph with no edges. Caro and Hansberg \cite{CaHa17} proved that $\iota(G, K_2) \leq n/3$ unless $G \simeq K_2$ or $G \simeq C_5$ (this was independently proved by \.{Z}yli\'{n}ski \cite{Z}, and the extremal graphs were investigated in \cite{LMS} and subsequently determined in \cite{BG}). Their seminal paper \cite{CaHa17} also featured several problems, some of which have been settled in the papers \cite{Borg1,Borgrsc,BFK} by the present second author and others (see also \cite{BFK2}). Domination and isolation have been particularly investigated for maximal outerplanar graphs \cite{BK, BK2, CaWa13, CaHa17, Ch75, DoHaJo16, DoHaJo17, HeKa18, LeZuZy17, Li16, MaTa96, To13, KaJi}, mostly due to connections with Chv\'{a}tal's Art Gallery Theorem \cite{Ch75}.

With a sharp upper bound established for each of $\iota(G,K_1)$ and $\iota(G,K_2)$, the next step was to consider isolation of each of the two connected $3$-vertex graphs $K_3$ and $P_3$. Again consider a connected $n$-vertex graph $G$. The sharp upper bound $n/4$ on $\iota(G,K_3)$ with $G \not\simeq K_3$ was established in \cite{Borg1} as a bound on the cycle isolation number, and is also given in \cite{BFK} as the special case $k = 3$ of the sharp upper bound $n/(k+1)$ on $\iota(G,K_k)$ with $G \not\simeq K_k$, in \cite{BorgIsolConnected2023} as a case of a result on isolation of connected $3$-edge graphs, and in \cite{Borgrsc} as a case of a result on isolation of regular graphs, stars and $k$-chromatic graphs. The significance of a $P_3$-isolating set $D$ is that the edge set of $G-N[D]$ is a matching, that is, no two edges of $G-N[D]$ have a common vertex (so $\Delta(G-N[D]) \leq 1$). The sharp upper bound $2n/7$ on $\iota(G,P_3)$ with $G$ not being a $\{P_3, C_3, C_6\}$-graph was established by Zhang and Wu in \cite{ZW}, and independently by the present second author in \cite{BorgIsolConnected2023}, which actually established the stronger inequality $\iota(G,P_3) \leq (4n-r)/14$ with $r$ being the number of leaves of $G$, and $G$ not being a copy of one of six particular graphs. The graphs attaining the bound are determined in each of \cite{CLWX} and \cite{CZZ}. Infinitely many of them have induced $6$-cycles. Zhang and Wu \cite{ZW} also showed that $\iota(G,P_3) \leq n/4$ if $G$ is not a $\{P_3,C_7,C_{11}\}$-graph and contains no cycles of order at most $6$. Huang, Zhang and Jin \cite{HZJ} improved this by showing that the same bound holds if $G$ is not a $\{P_3,C_3,C_7,C_{11}\}$-graph and $G$ either contains no $6$-cycles or contains $6$-cycles but not induced $5$-cycles or induced $6$-cycles. They asked if the bound still holds asymptotically for connected graphs having no induced $6$-cycles; see \cite[Problem~5.1]{HZJ}. Thus, the problem essentially is whether induced $6$-cycles solely account for the jump from the bound $n/4$ to the bound $2n/7$ on $\iota(G,P_3)$. In this paper, we solve this problem for subcubic graphs, which need to be treated differently from other graphs. We show that if $G$ is subcubic and has no induced $6$-cycles, then $\iota(G,P_3) \leq n/4$ unless $G$ is a copy of one of $12$ particular graphs. The bound is sharp. 

In order to state our result in full, we first recall an explicit construction in \cite{Borg1} that will be used for establishing an extremal case for any $n \neq 3$, and also construct the $12$ exceptions to the bound.

\begin{construction}[\cite{Borg1}] \label{Bconstruction} \emph{For any $n, k \in \mathbb{N}$ and any $k$-vertex graph $F$, we construct a connected $n$-vertex graph $B_{n,F}$ as follows.  If $n \leq k$, then let $B_{n,F} = P_n$. If $n \geq k+1$, then let $a_{n,k} = \big\lfloor\frac{n}{k+1}\big\rfloor$, let $b_{n,k} = n - ka_{n,k}$ (so $a_{n,k} \leq b_{n,k} \leq a_{n,k} + k$), let $F_1, \dots, F_{a_{n,k}}$ be copies of $F$ such that $P_{b_{n,k}}, F_1, \dots, F_{a_{n,k}}$ are pairwise vertex-disjoint, and let $B_{n,F}$ be the graph with $V(B_{n,F}) = [b_{n,k}] \cup \bigcup_{i=1}^{a_{n,k}} V(F_i)$ and $E(B_{n,F}) = E(P_{a_{n,k}}) \cup \{\{a_{n,k}, j\} \colon j \in [b_{n,k}] \setminus [a_{n,k}]\} \cup \bigcup_{i=1}^{a_{n,k}} (E(F_i) \cup \{\{i,v\} \colon v \in V(F_i)\})$.}
\end{construction}

Let $G_{7,1}, \dots, G_{7,6}$ be the $7$-vertex graphs whose drawings are given in Figure~1, let $G_{11}$ be the $11$-vertex graph whose drawing is given in Figure~2, and let $G_{15}$ be the $15$-vertex graph whose drawing is given in Figure~3. Let $\mathcal{G}_3 = \{P_3, C_3\}$, $\mathcal{G}_7 = \{C_7, G_{7,1}, \dots, G_{7,6}\}$, $\mathcal{G}_{11} = \{C_{11}, G_{11}\}$ and $\mathcal{G}_{15} = \{G_{15}\}$. Let $\mathcal{E} = \mathcal{G}_3 \cup \mathcal{G}_7 \cup \mathcal{G}_{11} \cup \mathcal{G}_{15}$. It is easy to see that for each $i \in \{3, 7, 11, 15\}$, $\iota(G,P_3) =(|V(G)|+1)/4$ for each $G \in \mathcal{G}_i$. Thus, if $G$ is an $\mathcal{E}$-graph, then $\iota(G,P_3) > |V(G)|/4$. We prove that $\mathcal{E}$-graphs are the only connected subcubic graphs for which this inequality holds. More precisely, in Section~\ref{section2}, we prove the following result.

\begin{theorem}\label{ThmP3}
If $G$ is a connected subcubic $n$-vertex graph that has no induced $6$-cycles and that is not an $\mathcal{E}$-graph, then 
$$\iota(G,P_3) \leq \floor*{\frac{n}{4}}.$$ 
Moreover, equality holds if $G = B_{n,P_3}$.
\end{theorem}

\noindent
The proof makes use of some of the latest ideas and methods in the study of isolation, such as divisibility considerations introduced in \cite{Borgisdom}. Our direct determination of the $\mathcal{G}_7$-graphs, or rather that these are the connected subcubic $7$-vertex graphs violating the bound in Theorem~\ref{ThmP3}, is corroborated by the results in \cite{CLWX, CZZ}. They are also confirmed by a computer check similar to that in \cite{BBS}, as are the $\mathcal{G}_{11}$-graphs. However, determining $G_{15}$ by hand has been inevitable and perhaps the most arduous part of this work as the enormity of the set of $15$-vertex graphs renders the use of the present computational power infeasible (the online data set in \cite{McKay}, used in \cite{BBS}, is for graphs of order at most $11$). Moreover, it was a surprise to discover that such a graph exists, especially when, unlike the cases of $\mathcal{G}_3$-graphs, $\mathcal{G}_7$-graphs and $\mathcal{G}_{11}$-graphs, a cycle with the same number of vertices does not violate the bound in Theorem~\ref{ThmP3}.

\begin{figure}[htb!]
\ffigbox[\textwidth]
    {
    \begin{floatrow}
    \ffigbox[\linewidth]
      {\captionof{subfigure}{$G_{7,1}$}
      \label{subfig:G1}}
      {\begin{tikzpicture} [scale =1, every edge/.style={draw=black, very thick}]
        \node[shape=circle,fill,draw=black, inner sep=0pt, scale=0.5, label={south west:$7$}] (7) at (0.5,0) {1};
            \node[shape=circle,fill,draw=black, inner sep=0pt, scale=0.5, label={north:$2$}] (2) at (1,1) {2};
            \node[shape=circle,fill,draw=black, inner sep=0pt, scale=0.5, label={north:$3$}] (3) at (2,1) {3};
            \node[shape=circle,fill,draw=black, inner sep=0pt, scale=0.5, label={north east:$4$}] (4) at (3,1) {4};
            \node[shape=circle,fill,draw=black, inner sep=0pt, scale=0.5, label={south east:$5$}] (5) at (2.5,0) {5};
            \node[shape=circle,fill,draw=black, inner sep=0pt, scale=0.5, label={south:$6$}] (6) at (1.5,-1) {6} ;
            \node[shape=circle,fill,draw=black, inner sep=0pt, scale=0.5, label={north west:$1$}] (1) at (0,1) {7};
            
            \draw (7) edge (2);
        		\draw (7) edge (6);
        		\draw (1) edge (7);
        		\draw (2) edge (3);
        		\draw (3) edge (4);
        		\draw (3) edge (5); 
       		\draw (4) edge (5);
        		\draw (5) edge (6);  
        		\draw (6) edge (7);
        
    \end{tikzpicture}}
    \ffigbox[\linewidth]
      {\captionof{subfigure}{$G_{7,2}$}
      \label{subfig:G2}}
      {\begin{tikzpicture} [scale =1, every edge/.style={draw=black, very thick}]
            \node[shape=circle,fill,draw=black, inner sep=0pt, scale=0.5, label={south west:$7$}] (7) at (0.5,0) {1};
            \node[shape=circle,fill,draw=black, inner sep=0pt, scale=0.5, label={north:$2$}] (2) at (1,1) {2};
            \node[shape=circle,fill,draw=black, inner sep=0pt, scale=0.5, label={north:$3$}] (3) at (2,1) {3};
            \node[shape=circle,fill,draw=black, inner sep=0pt, scale=0.5, label={north east:$4$}] (4) at (3,1) {4};
            \node[shape=circle,fill,draw=black, inner sep=0pt, scale=0.5, label={south east:$5$}] (5) at (2.5,0) {5};
            \node[shape=circle,fill,draw=black, inner sep=0pt, scale=0.5, label={south:$6$}] (6) at (1.5,-1) {6} ;
            \node[shape=circle,fill,draw=black, inner sep=0pt, scale=0.5, label={north west:$1$}] (1) at (0,1) {7};
            
            \draw (7) edge (2);
        		\draw (7) edge (6);
        		\draw (1) edge (7);
        		\draw (2) edge (3);
        		\draw (2) edge (1);
        		\draw (3) edge (4);
        		\draw (3) edge (5); 
       		\draw (4) edge (5);
        		\draw (5) edge (6);  
        		\draw (6) edge (7);
        \end{tikzpicture}}
    \end{floatrow}

    \vspace*{20pt}

    \begin{floatrow}
    \ffigbox[\linewidth]
      {\captionof{subfigure}{$G_{7,3}$}
      \label{subfig:G3}}
      {\begin{tikzpicture}[scale = 1, every edge/.style={draw=black, very thick}]
        \node[shape=circle,fill,draw=black, inner sep=0pt, scale=0.5, label={south west:$7$}] (7) at (0.5,0) {1};
            \node[shape=circle,fill,draw=black, inner sep=0pt, scale=0.5, label={north:$2$}] (2) at (1,1) {2};
            \node[shape=circle,fill,draw=black, inner sep=0pt, scale=0.5, label={north:$3$}] (3) at (2,1) {3};
            \node[shape=circle,fill,draw=black, inner sep=0pt, scale=0.5, label={north east:$4$}] (4) at (3,1) {4};
            \node[shape=circle,fill,draw=black, inner sep=0pt, scale=0.5, label={south east:$5$}] (5) at (2.5,0) {5};
            \node[shape=circle,fill,draw=black, inner sep=0pt, scale=0.5, label={south:$6$}] (6) at (1.5,-1) {6} ;
            \node[shape=circle,fill,draw=black, inner sep=0pt, scale=0.5, label={north west:$1$}] (1) at (0,1) {7};
            
            \draw (7) edge (2);
        		\draw (7) edge (6);
        		\draw (1) edge (7);
        		\draw (2) edge (3);
        		\draw (3) edge (4);
        		\draw (3) edge (5); 
       		\draw (4) edge (5);
        		\draw (5) edge (6);  
        		\draw (6) edge (7);
                \draw (1) edge[bend left = 1.8cm] (4);
        \end{tikzpicture}}
    \ffigbox[\linewidth]
      {\captionof{subfigure}{$G_{7,4}$}\label{subfig:G6}}
      {\begin{tikzpicture}[scale =1, every edge/.style={draw=black, very thick}]
            \node[shape=circle,fill,draw=black, inner sep=0pt, scale=0.5, label={south west:$7$}] (7) at (0.5,0) {1};
            \node[shape=circle,fill,draw=black, inner sep=0pt, scale=0.5, label={north:$2$}] (2) at (1,1) {2};
            \node[shape=circle,fill,draw=black, inner sep=0pt, scale=0.5, label={north:$3$}] (3) at (2,1) {3};
            \node[shape=circle,fill,draw=black, inner sep=0pt, scale=0.5, label={north east:$4$}] (4) at (3,1) {4};
            \node[shape=circle,fill,draw=black, inner sep=0pt, scale=0.5, label={south east:$5$}] (5) at (2.5,0) {5};
            \node[shape=circle,fill,draw=black, inner sep=0pt, scale=0.5, label={south:$6$}] (6) at (1.5,-1) {6} ;
            \node[shape=circle,fill,draw=black, inner sep=0pt, scale=0.5, label={north west:$1$}] (1) at (0,1) {7};
            
            \draw (7) edge (2);
        		\draw (7) edge (6);
        		\draw (1) edge (7);
        		\draw (2) edge (3);
        		\draw (2) edge (1);
        		\draw (3) edge (4);
        		\draw (3) edge (5); 
       		\draw (4) edge (5);
        		\draw (5) edge (6);  
        		\draw (6) edge (7);
                \draw (1) edge[bend left = 1.8cm] (4);
        \end{tikzpicture}}
    \end{floatrow}
    
    \vspace*{20pt}
    
    \begin{floatrow}
    \ffigbox[\linewidth]
      {\captionof{subfigure}{$G_{7,5}$}\label{subfig:G5}}
        {\begin{tikzpicture}[scale =1, every edge/.style={draw=black, very thick}]
          \node[shape=circle,fill,draw=black, inner sep=0pt, scale=0.5, label={south:$6$}] (1) at (-2,0) {1};
            \node[shape=circle,fill,draw=black, inner sep=0pt, scale=0.5, label={north east:$4$}] (2) at (-1,1) {2};
            \node[shape=circle,fill,draw=black, inner sep=0pt, scale=0.5, label={north:$3$}] (3) at (-2,1) {3};
            \node[shape=circle,fill,draw=black, inner sep=0pt, scale=0.5, label={north:$2$}] (4) at (-3,1) {4};
            \node[shape=circle,fill,draw=black, inner sep=0pt, scale=0.5, label={west:$1$}] (5) at (-4,0.5) {5};
            \node[shape=circle,fill,draw=black, inner sep=0pt, scale=0.5, label={south:$7$}] (6) at (-3,0) {6} ;
            \node[shape=circle,fill,draw=black, inner sep=0pt, scale=0.5, label={south east:$5$}] (7) at (-1,0) {7};
            
            \draw (1) edge (2);
        		\draw (1) edge (6);
        		\draw (1) edge (7);
        		\draw (2) edge (3);
        		\draw (3) edge (7);
        		\draw (3) edge (4); 
       		\draw (4) edge (5);
       		\draw (4) edge (6);
        		\draw (5) edge (6);  
        		\draw (6) edge (1);
       
        \end{tikzpicture}}
        
    \ffigbox[\linewidth]
      {\captionof{subfigure}{$G_{7,6}$}
      \label{subfig:G4}}
      {\begin{tikzpicture}[every edge/.style={draw=black, very thick}]
         \node[shape=circle,fill,draw=black, inner sep=0pt, scale=0.5, label={south:$6$}] (1) at (-2,0) {1};
            \node[shape=circle,fill,draw=black, inner sep=0pt, scale=0.5, label={north east:$4$}] (2) at (-1,1) {2};
            \node[shape=circle,fill,draw=black, inner sep=0pt, scale=0.5, label={north:$3$}] (3) at (-2,1) {3};
            \node[shape=circle,fill,draw=black, inner sep=0pt, scale=0.5, label={north:$2$}] (4) at (-3,1) {4};
            \node[shape=circle,fill,draw=black, inner sep=0pt, scale=0.5, label={west:$1$}] (5) at (-4,0.5) {5};
            \node[shape=circle,fill,draw=black, inner sep=0pt, scale=0.5, label={south:$7$}] (6) at (-3,0) {6} ;
            \node[shape=circle,fill,draw=black, inner sep=0pt, scale=0.5, label={south east:$5$}] (7) at (-1,0) {7};
            
                \draw (1) edge (2);
        		\draw (1) edge (6);
        		\draw (1) edge (7);
        		\draw (2) edge (3);
        		\draw (2) edge (7);
        		\draw (3) edge (7);
        		\draw (3) edge (4); 
       		\draw (4) edge (5);
       		\draw (4) edge (6);
        		\draw (5) edge (6);  
        		\draw (6) edge (1);
        \end{tikzpicture}}
    \end{floatrow}
     
    }{\caption{Subcubic $7$-vertex graphs with no induced $6$-cycles and with $P_3$-isolation number~$2$.}\label{FigG7}}
\end{figure}

\begin{figure}[htb!]
      {\caption{$G_{11}$, a subcubic $11$-vertex graph with no induced 6-cycles and with $P_3$-isolation number $3$.}
      \label{FigG11}}
      {\begin{tikzpicture}[scale = 1, every edge/.style={draw=black, very thick}]
           \node[shape=circle,fill,draw=black, inner sep=0pt, scale=0.5, label={south:$7$}] (7) at (0,0) {7};
        \node[shape=circle,fill,draw=black, inner sep=0pt, scale=0.5, label={north:$6$}] (6) at (0,1) {6};
        \node[shape=circle,fill,draw=black, inner sep=0pt, scale=0.5, label={north:$5$}] (5) at (-1,1) {5};
        \node[shape=circle,fill,draw=black, inner sep=0pt, scale=0.5, label={south:$8$}] (8) at (-1,0) {8};
        \node[shape=circle,fill,draw=black, inner sep=0pt, scale=0.5, label={south:$9$}] (9) at (-2,0) {9};
        \node[shape=circle,fill,draw=black, inner sep=0pt, scale=0.5, label={north:$4$}] (4) at (-2,1) {4} ;
        \node[shape=circle,fill,draw=black, inner sep=0pt, scale=0.375, label={south:$10$}] (10) at (-3,0) {10};
        \node[shape=circle,fill,draw=black, inner sep=0pt, scale=0.5, label={north:$3$}] (3) at (-3,1) {3};
        \node[shape=circle,fill,draw=black, inner sep=0pt, scale=0.375, label={south:$11$}] (11) at (-4,0) {11};
        \node[shape=circle,fill,draw=black, inner sep=0pt, scale=0.5, label={north:$2$}] (2) at (-4,1) {2};
        \node[shape=circle,fill,draw=black, inner sep=0pt, scale=0.5, label={west:$1$}] (1) at (-5,0.5) {1};
        
        \draw (1) edge (2);
        \draw (1) edge (11);
        \draw (2) edge (11);
        \draw (2) edge (3);
        \draw (3) edge (4);
        \draw (3) edge (9); 
        \draw (4) edge (5);
        \draw (4) edge (10);
        \draw (5) edge (6);  
        \draw (6) edge (7);
        \draw (7) edge (8);
        \draw (8) edge (9);
        \draw (9) edge (10);
        \draw (10) edge (11); 
        \end{tikzpicture}}
 {}
\end{figure}

\begin{figure}[htb!]

      {\caption{$G_{15}$, a subcubic $15$-vertex graph with no induced $6$-cycles and with $P_3$-isolation number $4$.}
      \label{FigG15}}
      {\begin{tikzpicture}[scale = 1, every edge/.style={draw=black, very thick}]
           \node[shape=circle,fill,draw=black, inner sep=0pt, scale=0.5, label={north:$8$}] (7) at (0,1) {7};
        \node[shape=circle,fill,draw=black, inner sep=0pt, scale=0.5, label={north:$7$}] (6) at (-1,1) {6};
        \node[shape=circle,fill,draw=black, inner sep=0pt, scale=0.5, label={north:$6$}] (5) at (-2,1) {5};
        \node[shape=circle,fill,draw=black, inner sep=0pt, scale=0.5, label={east:$9$}] (8) at (1,0.5) {8};
        \node[shape=circle,fill,draw=black, inner sep=0pt, scale=0.5, label={south:$10$}] (9) at (0,0) {9};
        \node[shape=circle,fill,draw=black, inner sep=0pt, scale=0.5, label={north:$4$}] (4) at (-4,1) {4} ;
        \node[shape=circle,fill,draw=black, inner sep=0pt, scale=0.375, label={south:$11$}] (10) at (-1,0) {10};
        \node[shape=circle,fill,draw=black, inner sep=0pt, scale=0.5, label={north:$3$}] (3) at (-5,1) {3};
        \node[shape=circle,fill,draw=black, inner sep=0pt, scale=0.375, label={south:$12$}] (11) at (-2,0) {11};
        \node[shape=circle,fill,draw=black, inner sep=0pt, scale=0.5, label={north:$2$}] (2) at (-6,1) {2};
        \node[shape=circle,fill,draw=black, inner sep=0pt, scale=0.5, label={west:$1$}] (1) at (-7,0.5) {1};
        \node[shape=circle,fill,draw=black, inner sep=0pt, scale=0.375, label={north:$5$}] (12) at (-3,1) {12};
        \node[shape=circle,fill,draw=black, inner sep=0pt, scale=0.375, label={south:$13$}] (13) at (-4,0) {13};
        \node[shape=circle,fill,draw=black, inner sep=0pt, scale=0.375, label={south:$14$}] (14) at (-5,0) {14};
        \node[shape=circle,fill,draw=black, inner sep=0pt, scale=0.375, label={south:$15$}] (15) at (-6,0) {15}; 
        
        \draw (1) edge (2);
        \draw (1) edge (15);
        \draw (2) edge (15);
        \draw (2) edge (3);
        \draw (3) edge (4);
        \draw (3) edge (13); 
        \draw (4) edge (12);
        \draw (12) edge (5);
        \draw (4) edge (14);
        \draw (5) edge (6);
        \draw (5) edge (10);  
        \draw (6) edge (7);
        \draw (6) edge (11);
        \draw (7) edge (8);
        \draw (7) edge (9);
        \draw (8) edge (9);
        \draw (9) edge (10);
        \draw (10) edge (11);
        \draw (11) edge (13);
        \draw (13) edge (14);
        \draw (14) edge (15);
        \end{tikzpicture}}
    
 {}
\end{figure}

\section{Proof of Theorem~\ref{ThmP3}} \label{section2}

We now prove Theorem~\ref{ThmP3}. We often use the following two lemmas. 

\begin{lemma}[\cite{Borg1}] \label{LemIsol1} If $G$ is a graph, $\mathcal{F}$ is a set of graphs, $X \subseteq V(G)$ and $Y \subseteq N[X]$, then $$\iota(G,\mathcal{F}) \leq |X| + \iota(G-Y,\mathcal{F}).$$
\end{lemma}

A \emph{component} of a graph $G$ is a maximal connected subgraph of $G$. Clearly, the components of $G$ are pairwise vertex-disjoint, and their union is $G$. Let ${\rm C}(G)$ denote the set of components of $G$. Let ${\rm C}_3(G)$ denote the set of components of $G$ containing a $3$-path.
    
\begin{lemma}[\cite{Borg1, Borgrsc}] \label{LemIsol2} For any graph $G$ and any set $\mathcal{F}$ of connected graphs, $$\iota(G,\mathcal{F}) = \sum_{H \in {\rm C}(G)} \iota(H,\mathcal{F}).$$
\end{lemma}

We denote the set of $P_3$-isolating sets of a graph $G$ by ${\rm Is}(G, P_3)$.

\begin{lemma}\label{LemP3G5} If a graph $G$ has at most $5$ vertices, then $\iota(G,P_3)\leq 1$.
\end{lemma}
\proof{If $G$ contains no $P_3$-copy, then $\iota(G,P_3) = 0$. Suppose that $G$ contains a $P_3$-copy. Then, $G$ has a vertex $v$ such that $d(v) \geq 2$. Thus, $|V(G-N[v])| \leq 2$, and hence $\{v\} \in {\rm Is}(G, P_3)$. \qed}

\begin{lemma}\label{LemDelta2}  If $G$ is a connected $n$-vertex graph that is not a $\{P_3,C_3,C_6,C_7,C_{11}\}$-graph, and $\Delta(G) \leq 2$, then $\iota(G,P_3) \leq n/4$.
\end{lemma}
\proof{If $n \leq 2$, then $\iota(G,P_3) = 0$. If $3 \leq n \leq 5$, then the result is given by Lemma~\ref{LemP3G5}. Suppose $n \geq 6$. Since $G$ is connected and $\Delta(G) \leq 2$, we have $G \simeq P_n$ or $G\simeq C_n$. If $G = P_n$, then $\{4k \colon k \in [\floor*{n/4}]\} \in {\rm Is}(G, P_3)$. If $G = C_n$, then $\{5k-4 \colon k \in \mathbb{N}, \, 5k-4 \leq n\}$ is a $P_3$-isolating set of $G$ of size $\floor*{(n+4)/5}$. Clearly, $\floor*{(n+4)/5} \leq n/4$ if $n \in \{8, 9, 10, 12, 13, 14, 15\}$, and $(n+4)/5 \leq n/4$ if $n \geq 16$.  \qed}
\\

The next lemma is Theorem~\ref{ThmP3} for $n \leq 15$. In order to bring the key ideas to the fore, we first prove Theorem~\ref{ThmP3} by assuming this lemma, and then we provide the proof of the lemma, which is particularly laborious for the case $n = 15$. 

\begin{lemma}\label{ThmP3small}
If $G$ is a connected subcubic $n$-vertex graph with $n \leq 15$, and $G$ has no induced $6$-cycles and is not an $\mathcal{E}$-graph, then $\iota(G,P_3) \leq n/4$.
\end{lemma}

The following observations are fairly straightforward.

\begin{observation}\label{ObsE} Let $G \in \mathcal{E}$. \\ 
(a) $\iota(G,P_3) = (|V(G)|+1)/4$.\\
(b) If $v \in V(G)$ and $(G,v) \notin \{(P_3,2), (G_{7,1},7)\}$, then $G-v$ is connected.
\end{observation}

\begin{observation}\label{ObsG7} Let $G \in \mathcal{G}_7$ such that $G \neq C_7$.\\
(a) For each $v \in V(G)$ with $d(v) \leq 2$, there exists some $v' \in V(G) \setminus N[v]$ such that $d(v') = 3$ and $G-N[v']$ is connected. \\
(b) No two vertices of $G$ of degree $1$ or $2$ are adjacent.  
\end{observation}

\begin{observation}\label{ObsG11} Let $G \in \mathcal{G}_{11} \cup \mathcal{G}_{15}$. \\
(a) $\delta(G) = 2$. \\
(b) For each $v \in V(G)$ with $d(v) \leq 2$, $G-N[v]$ is connected. \\
(c) If $G \in \mathcal{G}_{15}$, then $\min\{d(v,w) \colon v, w \in V(G), \, v \neq w, \, d(v) = d(w) = 2\} = 4$.
\end{observation}

\begin{observation}\label{ObsG7i} Let $G \in \{G_{7,i} \colon i \in \{1, 2, 3, 5\}\}$, and let $v \in V(G)$ with $d(v) \leq 2$. There exists some $v' \in V(G)\setminus N[v]$ such that $d(v') = 3$, $G-N[v']$ is a $\mathcal{G}_3$-graph, and one of the following statements holds: \\
(a) $G-N[v'] \simeq P_3$ and there exist $y,y' \in V(G-N[v'])$ such that $d(y) = 2$, $d_{G-N[v']}(y) = 1$, $d(y') = 3$ and $d_{G-N[v']}(y') = 2$.\\
(b) $G-N[v'] \simeq C_3$ and $d(y) = d(y') = 3$ for some $y,y' \in V(G-N[v'])$ with $y \neq y'$.
\end{observation}

\begin{lemma}\label{LemG-v} If $G \in \mathcal{E}$ and $v \in V(G)$ with $d(v) \leq 2$, then $\iota(G-v,P_3) \leq (|V(G)|-3)/4$, and if $G \notin \mathcal{G}_3$, then $G-v$ is connected.
\end{lemma}
\proof{The result is trivial if $G \in \mathcal{G}_3$. Suppose $G \notin \mathcal{G}_3$. We have $|V(G)| = 4k + 3$ for some $k \in [3]$. Since $|V(G-v)| = 4k+2$, $G-v$ is not an $\mathcal{E}$-graph. Since $d_{G_{7,1}}(7) = 3$, $(G, v) \neq (G_{7,1}, 7)$. By Observation~\ref{ObsE}, $G-v$ is connected. By Lemma~\ref{ThmP3small}, $\iota(G-v,P_3) \leq \floor{|V(G-v)|/4} = (|V(G)|-3)/4$. \qed}
\\

If $v \in V(G)$, $x \in N(v)$, $H \in {\rm C}(G - N[v])$, and $xy_{x,H} \in E(G)$ for some $y_{x,H} \in V(H)$, then, as in \cite{Borg1}, we say that $x$ \emph{is linked to} $H$ and that $H$ \emph{is linked to} $x$.
\\

\noindent
\textbf{Proof of Theorem~\ref{ThmP3}.} Let $B$ be the graph $B_{n,P_3}$ defined in Construction~\ref{Bconstruction}. If $n \leq 2$, then $\iota(B,P_3) = 0$. If $n = 3$, then $B \simeq P_3$. Suppose $n \geq 4$. Since $[a_{n,3}] \in {\rm Is}(B, P_3)$, $\iota(B,P_3) \leq a_{n,3}$. If $D \in {\rm Is}(B, P_3)$, then for each $i \in [a_{n,3}]$, $D \cap (V(F_i)\cup\{i\}) \neq \emptyset$ as $B-N_B[D]$ contains no $P_3$-copy. Thus, $\iota(B,P_3) = a_{n,3} = \floor{n/4}$. Therefore, the bound is attained if $G = B$.

We now use induction on $n$ to prove the inequality in the theorem. Since $\iota(G,P_3)$ is an integer, it suffices to prove that $\iota(G,P_3) \leq n/4$. This holds by Lemma~\ref{ThmP3small} if $n \leq 15$, and by Lemma~\ref{LemDelta2} if $\Delta(G) \leq 2$. Suppose $n \geq 16$ and $\Delta(G) = 3$. Let $v \in V(G)$ with $d(v) = 3$. Since $n \geq 16$, $N[v] \neq V(G)$. Let $G' = G-N[v]$. Then, $V(G') \neq \emptyset$. Let $\mathcal{H} = {\rm C}(G')$, $\mathcal{H}' = \{H \in \mathcal{H} \colon \iota(H,P_3) > |V(H)|/4\}$ and $\mathcal{H}^* = \mathcal{H} \setminus \mathcal{H}'$. By the induction hypothesis, each member of $\mathcal{H}'$ is an $\mathcal{E}$-graph. 
    
If $\mathcal{H}' = \emptyset$, then by Lemmas~\ref{LemIsol1} and~\ref{LemIsol2}, 
\begin{align*}
\iota(G,P_3) &\leq 1 + \iota(G',P_3) \leq \frac{|N[v]|}{4} + \sum_{H\in\mathcal{H}} {\frac{|V(H)|}{4}}= \frac{n}{4}.
\end{align*}
Suppose $\mathcal{H}' \neq \emptyset$. Since $G$ is connected, each member of $\mathcal{H}$ is linked to at least one member of $N(v)$. For each $x \in N(v)$, let $\mathcal{H}_{x} = \{H \in \mathcal{H} \colon H$ is linked to $ x\}$, $\mathcal{H}'_{x} = \{H \in \mathcal{H}' \colon H $ is linked to $x\}$ and $\mathcal{H}_{x}^* = \{H \in \mathcal{H}^* \colon H$ is linked to $x$ only$\}$. For each $H \in \mathcal{H}^*$, let $D_H$ be a $P_3$-isolating set of $H$ of size $\iota(H,P_3)$. 

For any $x \in N(v)$ and any $H \in \mathcal{H}_x'$, let $H_x' = H - y_{x,H}$. Since $\Delta(G) = 3$, $d_H(y_{x,H}) \leq 2$.
Let $D_{x,H}$ be a smallest $P_3$-isolating set of $H_x'$. By Lemma \ref{LemG-v}, 
\begin{equation} |D_{x,H}| \leq \frac{|V(H)|-3}{4}. \label{DxHineq}
\end{equation}

\noindent
\textbf{Case 1:} \emph{$|\mathcal{H}'_x| \geq 2$ for some $x \in N(v)$.} (Since $\Delta(G) = 3$ and $v \in N(x)$, we actually have $|\mathcal{H}'_x| = 2$, but this fact is not used here.) For each $H \in \mathcal{H}' \setminus \mathcal{H}'_x$, let $x_H \in N(v)$ such that $H$ is linked to $x_H$. Let $X = \{x_H \colon H \in \mathcal{H}' \setminus \mathcal{H}'_x\}$. Note that $x \notin X$. Let 
\[D = \{v, x\} \cup X \cup \left( \bigcup_{H \in \mathcal{H}_x'} D_{x,H} \right) \cup \left( \bigcup_{H \in \mathcal{H}' \setminus \mathcal{H}_x'} D_{x_H,H} \right) \cup \left( \bigcup_{H \in \mathcal{H}^*} D_{H} \right).\]
We have $V(G) = N[v] \cup \bigcup_{H \in \mathcal{H}} V(H)$, $y_{x,H} \in N[x]$ for each $H \in \mathcal{H}'_x$, and $y_{x_H,H} \in N[x_H]$ for each $H \in \mathcal{H}' \setminus \mathcal{H}'_x$, so $D \in {\rm Is}(G, P_3)$. We have
\begin{align*} \iota(G, P_3) &\leq |D| = 2 + |X| + \sum_{H \in \mathcal{H}_x'} |D_{x,H}| + \sum_{H \in \mathcal{H}' \setminus \mathcal{H}_x'} |D_{x_H,H}| + \sum_{H \in \mathcal{H}^*} |D_{H}| \nonumber \\
&\leq \frac{8 + 4|X|}{4} + \sum_{H \in \mathcal{H}_x'} \frac{|V(H)| - 3}{4} + \sum_{H \in \mathcal{H}' \setminus \mathcal{H}_{x}'} \frac{|V(H)| - 3}{4} + \sum_{H \in \mathcal{H}^*} \frac{|V(H)|}{4} \nonumber \\
&= \frac{(2 + |X|) + 3(2 - |\mathcal{H}_x'|) + 3(|X| - |\mathcal{H}' \setminus \mathcal{H}'_x|)}{4} + \sum_{H \in \mathcal{H}} \frac{|V(H)|}{4}.
\end{align*}
Thus, since $2 \leq |\mathcal{H}'_x|$ and $|X| \leq |\mathcal{H}' \setminus \mathcal{H}'_x|$, $\iota(G, P_3) \leq \frac{2 + |X|}{4} + \sum_{H \in \mathcal{H}} \frac{|V(H)|}{4} \leq \frac{n}{4}$.\medskip

\noindent
\textbf{Case 2:} \emph{$|\mathcal{H}'_x| \leq 1$ for each $x \in N(v)$.} For each $H \in \mathcal{H}'$, let $x_H \in N(v)$ such that $H$ is linked to $x_H$.\medskip 

\noindent
\textbf{Case 2.1:} \emph{Some member $H$ of $\mathcal{H}'$ is linked to $x_H$ only.} Let $G^* = G - (\{x_H\} \cup V(H))$. Then, $G^*$ has a component $G_v^*$ such that $N[v] \setminus \{x_H\} \subseteq V(G_v^*)$, and any other component of $G^*$ is a member of $\mathcal{H}_{x_H}^*$. Let $D^*$ be a smallest $P_3$-isolating set of $G_v^*$. Let $D = D^* \cup \{x_H\} \cup D_{x_H,H} \cup \bigcup_{I \in \mathcal{H}_{x_H}^*} D_I$. Since $D \in {\rm Is}(G, P_3)$,
\begin{align*} \iota(G, P_3) &\leq |D^*| + 1 + |D_{x_H,H}| + \sum_{I \in \mathcal{H}_{x_H}^*} |D_I| \\
&\leq \iota(G_v^*, P_3) + \frac{|\{x_H\} \cup V(H)|}{4} + \sum_{I \in \mathcal{H}_{x_H}^*} \frac{|V(I)|}{4}
\end{align*}
by (\ref{DxHineq}). This yields $\iota(G, P_3) \leq n/4$ if $\iota(G_v^*, P_3) \leq |V(G_v^*)|/4$. Suppose $\iota(G_v^*, P_3) > |V(G_v^*)|/4$. By the induction hypothesis, $G_v^*$ is an $\mathcal{E}$-graph. Since $v \in N[x_H]$, $d_{G_v^*}(v) = 2$. By Lemma \ref{LemG-v}, $G_v^* - v$ has a $P_3$-isolating set $D_v$ with $|D_v| \leq (|V(G_v^*)|-3)/4$. Since $v \in N[x_H]$, $(D \setminus D^*) \cup D_v \in {\rm Is}(G, P_3)$, and hence
\begin{equation*} \iota(G, P_3) \leq |D_v| + \frac{|\{x_H\} \cup V(H)|}{4} + \sum_{I \in \mathcal{H}_{x_H}^*} \frac{|V(I)|}{4} < \frac{|V(G_v^*)|}{4} + \frac{n - |V(G_v^*)|}{4} = \frac{n}{4}.
\end{equation*}
\noindent
\textbf{Case 2.2:} \emph{For each $H \in \mathcal{H}'$, $H$ is linked to some $x_H' \in N(v) \setminus \{x_H\}$.} Since $|\mathcal{H}_x'| \leq 1$ for each $x \in N(v)$, no two members of $\mathcal{H}'$ are linked to the same neighbour of $v$. Thus, since $d(v) = \Delta(G) = 3$ and $\mathcal{H}' \neq \emptyset$, $|\mathcal{H}'| = 1$. Let $H_1$ be the member of $\mathcal{H}'$. Then, $H_1$ is an $\mathcal{E}$-graph. Let $x_1 = x_{H_1}$ and $x_1' = x_{H_1}'$. Let $w$ be the unique member of $N(v) \setminus \{x_1,x_1'\}$. Let $Y = N[v] \cup V(H_1)$. Then, $V(G) = Y \cup \bigcup_{H \in \mathcal{H}^*} V(H)$.\medskip

\noindent
\textbf{Case 2.2.1:} \emph{$H_1$ is a $\mathcal{G}_{15}$-graph.} We may assume that $H_1 = G_{15}$. Thus, $n \geq 19$. Let $y_i = i$ for each $i\in [15]$. Since $\Delta(G) = 3$, we have $d_G(y_i) = d_{H_1}(y_i)$ for each $i \in [15] \setminus \{1, 5, 9\}$, and $|N(y_j) \cap N(v)| \leq 1$ for each $j \in \{1, 5, 9\}$. It follows that $y_{x_1,H_1} \neq y_{x_1',H_1}$ and that $y_{x_1,H_1} \in \{y_1, y_9\}$ or $y_{x_1',H_1} \in \{y_1, y_9\}$ (otherwise, $y_{x_1, H_1} = y_{x_1', H_1} = y_5$, a contradiction). By symmetry, we may assume that $y_{x_1,H_1} = y_1$. Thus, $y_{x_1',H_1} \in \{y_5,y_9\}$. Let $G^* = G-N[y_{13}]$. Then, $G^*$ is connected, $|V(G^*)| \geq 15$ and $d_{G^*}(y_4) = 1$. By Observation \ref{ObsG11}~(a), $G^*$ is not a $\mathcal{G}_{15}$-graph. By Lemma~\ref{LemIsol1} and the induction hypothesis, $\iota(G,P_3) \leq 1 + \iota(G^*,P_3) \leq 1 + (n-4)/4 = n/4$.\medskip

\noindent
\textbf{Case 2.2.2:} \emph{$H_1$ is a $\mathcal{G}_{11}$-graph.} Let $y_i = i$ for each $i \in [11]$. 

Suppose $H_1 \simeq G_{11}$. We may assume that $H_1 = G_{11}$. Similarly to Case~2.2.1, at least one of $y_5, y_6, y_7$ and $y_8$ is adjacent to some $x \in \{x_1, x_1'\}$ (otherwise, $y_{x_1, H_1} = y_{x_1', H_1} = y_1$, a contradiction). Let $G^* = G - N[y_2]$. Then, $G^*$ is connected, $|V(G^*)| \geq 12$, $d_{G^*}(y_4) = d_{G^*}(y_{10}) = 2$ and $y_4y_{10} \in E(G^*)$. By Observation~\ref{ObsG11}~(c), $G^*$ is not a $\mathcal{G}_{15}$-graph. We obtain $\iota(G,P_3) \leq n/4$ as in Case~2.2.1.

Now suppose $H_1 \simeq C_{11}$. We may assume that $H_1 = C_{11}$ and $y_1 = y_{x_1,H_1}$. Since $|Y| = 15$ and $n \geq 16$, there exists some $z \in V(G) \setminus Y$ such that $uz \in E(G)$ for some $u \in N(v)$. Let $u' \in \{x_1, x_1'\} \setminus\{u\}$. We may assume that $u' = x_1$. Since $\Delta(G) = 3$, $N(y_1) = \{x_1, y_2, y_{11}\}$. Let $I = G[Y] - N[y_1]$. Then, $V(I) = \{v, x_1', w\} \cup (V(H_1)\setminus N[y_1])$. 

Suppose that $I$ is connected. Let $G^* = G - N[y_1]$. Then, $G^*$ has a component $G_v^*$ such that $V(I) \cup \{z\} \subseteq V(G_v^*)$, and ${\rm C}(G^*) = \{G_v^*\} \cup \mathcal{H}_{x_1}^*$. Note that $|V(G_v^*)| \geq 12$, $d_{G_v^*}(v) = 2$, $d_{G_v^*}(y_3) \leq 2$ and $d_{G_v^*}(y_{10}) \leq 2$. If $y_3 \notin N(x_1') \cup N(w)$, then $d_{G_v^*}(y_3) = 1$, and hence, by Observation~\ref{ObsG11}~(a), $G_v^*$ is not a $\mathcal{G}_{15}$-graph. If, on the other hand, $y_3 \in N(x_1') \cup N(w)$, then $d_{G_v^*}(v) = d_{G_v^*}(y_3) = 2$ and $d(v, y_3) = 2$, and hence, by Observation~\ref{ObsG11}~(c), $G_v^*$ is not a $\mathcal{G}_{15}$-graph. Thus, $G_v^*$ is not an $\mathcal{E}$-graph. By Lemma~\ref{LemIsol1}, Lemma~\ref{LemIsol2} and the induction hypothesis, 
\begin{align*} \iota(G, P_3) &\leq 1 + \iota(G^*, P_3) = 1 + \iota(G_v^*, P_3) + \sum_{H \in \mathcal{H}_{x_1}^*} \iota(H, P_3) \\
&\leq \frac{|N[y_1]|}{4} + \frac{|V(G_v^*)|}{4} + \sum_{H \in \mathcal{H}_{x_1}^*} \frac{|V(H)|}{4} = \frac{n}{4}.
\end{align*}

Now suppose that $I$ is not connected. Then, the members of ${\rm C}(I)$ are the $\mathcal{G}_3$-graph $G[\{v,x_1',w\}]$ and the $8$-vertex path $G[V(H_1)\setminus N[y_1]]$, and $y_{x_1',H_1} \in \{y_2, y_{11}\}$. Let $\mathcal{I} = \{H \in \mathcal{H}^* \colon H \text{ is only linked to one or both members of }\{x_1,x_1'\}\}$ and $G^* = G - (Y \setminus \{v,w\})$. Then, $G^*$ has a component $G_v^*$ such that $\{v,w\} \subseteq V(G_v^*)$, and ${\rm C}(G^*) = \{G_v^*\} \cup \mathcal{I}$. If $G_v^*$ is an $\mathcal{E}$-graph, then $\iota(G_v^*, P_3) = (|V(G_v^*)|+1)/4$ by Observation~\ref{ObsE}; otherwise, $\iota(G_v^*, P_3) \leq |V(G_v^*)|/4$ by the induction hypothesis. Let $D^*$ be a $P_3$-isolating set of $G_v^*$ of size $\iota(G_v^*, P_3)$. If $y_{x_1',H_1} = y_2$, then let $D' = \{y_1,y_2,y_7\}$. If $y_{x_1',H_1} = y_{11}$, then let $D' = \{y_1,y_6,y_{11}\}$. Then, $D' \cup D^* \cup \bigcup_{H \in \cI} D_H \in {\rm Is}(G, P_3)$. Therefore,
\begin{equation*} \iota(G, P_3) \leq \frac{|Y \setminus \{v,w\}| - 1}{4} + \frac{|V(G_v^*)|+1}{4} + \sum_{H \in \mathcal{I}} \frac{|V(H)|}{4} = \frac{n}{4}.
\end{equation*}

\noindent
\textbf{Case~2.2.3:} \emph{$H_1$ is a $\mathcal{G}_7$-graph.} Since $H_1$ is linked to $2$ distinct neighbours of $v$, $H_1$ has a leaf or at least two vertices of degree $2$ in $H_1$. Thus, $H_1$ is a $(\mathcal{G}_7 \setminus \{G_{7,4}, G_{7,6}\})$-graph. 

Suppose $H_1 \not\simeq C_7$. By Observation~\ref{ObsG7i}, there exists some $y^* \in V(H_1) \setminus N[y_{x_1,H_1}]$ such that $d_{H_1}(y^*) = 3$, $H_1 -N_{H_1}[y^*]$ is a $\mathcal{G}_3$-graph, and $H_1 -N_{H_1}[y^*]$ has two distinct vertices $z$ and $z'$ such that $d_{H_1 - N_{H_1}[y^*]}(z) = d_{H_1}(z) - 1 \leq 2$ and $d_{H_1 - N_{H_1}[y^*]}(z') = d_{H_1}(z') - 1 = 2$. Thus, $N_G(y^*) = N_{H_1}(y^*)$, $1 \leq d_{G-N_G[y^*]}(z) \leq \Delta(G)  - 1 = 2$, $N_G(z') = N_{H_1}(z')$ and $zz' \in E(H_1)$. Let $G^* = G - N[y^*]$. Then, $G^*$ has a leaf or two adjacent vertices of degree $2$ in $G^*$. By Observation~\ref{ObsG11}, $G^*$ is not a $\mathcal{G}_{15}$-graph. Since $|V(G^*)| = n-4 \geq 12$, $G^*$ is not an $\mathcal{E}$-graph. Since $x_1y_{x_1,H_1} \in E(G^*)$, $G^*$ is connected. By Lemma~\ref{LemIsol1} and the induction hypothesis, $\iota(G, P_3) \leq 1 + \iota(G^*,P_3) \leq 1 + (n-4)/4 = n/4$.

Now suppose $H_1 \simeq C_7$. We may assume that $H_1 = C_7$. Let $y_i = i$ for each $i\in [7]$. We may assume that $y_1 = y_{x_1,H_1}$. Since $|Y| = 11$ and $n \geq 16$, there exists some $z \in V(G)\setminus Y$ such that $uz \in E(G)$ for some $u \in N(v)$. Let $u' \in \{x_1,x_1'\}\setminus\{u\}$. We may assume that $u' = x_1$. Let $I = G[Y]-N[y_1]$, and let $J$ be the $4$-vertex path $H_1 - N_{H_1}[y_1]$. Then, $V(I) = \{v, x_1', w\} \cup V(J)$. Let $G^* = G-N[y_1]$. Since $\Delta(G) = 3$, $N(y_1) = \{x_1, y_2, y_7\}$.

Suppose that $I$ is connected. Then, $G^*$ has a component $G_v^*$ such that $V(I) \cup \{z\} \subseteq V(G_v^*)$, and ${\rm C}(G^*) = \{G_v^*\} \cup \mathcal{H}_{x_1}^*$. Note that $|V(G_v^*)| \geq 8$, $d_{G_v^*}(v)=2$, $d_{G_v^*}(y_3)\leq 2$ and $d_{G_v^*}(y_6)\leq 2$. If $y_3 \notin N(x_1') \cup N(w)$, then $d_{G_v^*}(y_3) = 1$, and hence, by Observation~\ref{ObsG11}~(a), $G_v^*$ is not a $(\mathcal{G}_{11} \cup \mathcal{G}_{15})$-graph. If, on the other hand, $y_3 \in N(x_1') \cup N(w)$, then $d_{G_v^*}(v) = d_{G_v^*}(y_3) = 2$, $d(v, y_3) = 2$, $G_v^* - N_{G_v^*}[v]$ is not connected (as $z \notin Y$), and hence, by parts (b) and (c) of Observation~\ref{ObsG11}, $G_v^*$ is not a $(\mathcal{G}_{11} \cup \mathcal{G}_{15})$-graph. Thus, $G_v^*$ is not an $\mathcal{E}$-graph. As in Case~2.2.2, we obtain $\iota(G,P_3) \leq n/4$.

Now suppose that $I$ is not connected. Then, $G^*$ has a component $G_v^*$ such that $v, x_1', w, z \in V(G_v^*)$, ${\rm C}(G^*) = \{G_v^*, J\} \cup \mathcal{H}_{x_1}^*$ and $y_{x_1',H_1} \in \{y_2, y_7\}$. Let $D^*$ be a $P_3$-isolating set of $G_v^*$ of size $\iota(G_v^*, P_3)$. Then, $D^* \cup \{y_5\} \cup \bigcup_{H \in \mathcal{H}_{x_1}^*} D_H \in {\rm Is}(G^*, P_3)$. If $G_v^*$ is not an $\mathcal{E}$-graph, then by Lemma~\ref{LemIsol1} and the induction hypothesis, 
\begin{align*} \iota(G, P_3) &\leq 1 + \iota(G^*, P_3) \leq 2 + |D^*| + \sum_{H \in \mathcal{H}_{x_1}^*} |D_H| \\
&\leq \frac{|V(H_1) \cup \{x_1\}|}{4} + \frac{|V(G_v^*)|}{4} + \sum_{H \in \mathcal{H}_{x_1}^*} \frac{|V(H)|}{4} = \frac{n}{4}.
\end{align*}
Suppose that $G_v^*$ is an $\mathcal{E}$-graph. Since $|V(G_v^*)| \geq 4$, $G_v^*$ is not a $\mathcal{G}_3$-graph. Since $d_{G_v^*}(v) = 2$ and $d_{G_v^*}(x_1') \leq 2$, we obtain from Observations~\ref{ObsG7}~(b) and~\ref{ObsG11} (c) that $G_v^*$ is a $\{C_7, C_{11}, G_{11}\}$-graph. If $G_v^*$ is a $\mathcal{G}_{11}$-graph, then we obtain $\iota(G, P_3) \leq n/4$ by the argument for Case~2.2.2 with $y_1$ and $G_v^*$ taking the roles of $v$ and $H_1$ there, respectively (note that we now have that $G_v^*$ is linked to the neighbours $x_1$ and $y_{x_1',H_1}$ of $y_1$). Now suppose $G_v^* \simeq C_7$. Thus, $G_v^*$ has $4$ distinct vertices $u_1, u_2, u_3, u_4 \notin \{v, x_1', w\}$ such that $E(G_v^*) = \{vx_1', x_1'u_1, u_1u_2, u_2u_3, u_3u_4, u_4w, wv\}$. Recall that $y_{x_1',H_1} \in \{y_2, y_7\}$. We may assume that $y_{x_1',H_1} = y_7$. Since $\Delta(G) = 3$, $N(x_1') = \{v, u_1, y_7\}$. Let $n^* = \sum_{H \in \mathcal{H}_{x_1}^*} |V(H)|$. Since ${\rm C}(G^*) = \{G_v^*, J\} \cup \mathcal{H}_{x_1}^*$, $n = |N[y_1]| + |V(G_v^*)| + |V(J)| + n^* = 15 + n^*$. Thus, $\mathcal{H}_{x_1}^* \neq \emptyset$ as $n \geq 16$. Let $H^* \in \mathcal{H}_{x_1}^*$. Since $\Delta(G) = 3$, $N(x_1) = \{v, y_1, y_{x_1,H^*}\}$ and $\mathcal{H}_{x_1}^* = \{H^*\}$. Let $G^{\dagger} = G - N[\{x_1',u_4\}]$. The components of $G^{\dagger}$ are $(\{u_2\}, \emptyset)$ and a graph $H^{\dagger}$ with $V(H^{\dagger}) = V(H_1 - y_7) \cup \{x_1\} \cup V(H^*)$. Since $|V(H^{\dagger})| \geq 8$ and $d_{H^{\dagger}}(y_6) = 1$, $H^{\dagger}$ is not an $\mathcal{E}$-graph, so $\iota(H^{\dagger}, P_3) \leq |V(H^{\dagger})|/4$ by the induction hypothesis. By Lemmas~\ref{LemIsol1} and~\ref{LemIsol2}, 
\begin{equation*} \iota(G, P_3) \leq 2 + \iota(G^{\dagger}, P_3) = 2 + \iota(H^{\dagger}, P_3) \leq \frac{|V(G_v^*) \cup \{y_7\}|}{4} + \frac{|V(H^{\dagger})|}{4} = \frac{n}{4}.
\end{equation*}

\noindent
\textbf{Case 2.2.4:} \emph{$H_1$ is a $\mathcal{G}_3$-graph.} Let $y_1 = y_{x_1,H_1}$ and $y_1' = y_{x_1',H_1}$. 

Suppose that 
\begin{equation} \mbox{no vertex of $H_1$ of degree $2$ in $H_1$ is adjacent to a neighbour of $v$ (in $G$).} \label{2.2.4-1}
\end{equation} 
Thus, $d_{H_1}(y_1) = d_{H_1}(y_1') = 1$ and $H_1\simeq P_3$. 

Suppose $y_1 = y_1'$. Then, $H_1 = (\{y_1, y', y^*\}, \{y_1y', y'y^*\})$ for some distinct $y'$ and $y^*$ in $V(H_1) \setminus \{y_1\}$, $N(y_1) = \{x_1, x_1', y'\}$, and $N(y') = \{y_1, y^*\}$ by (\ref{2.2.4-1}). Let $G^* = G-N[y_1]$ and $\mathcal{I} = \{H \in \mathcal{H}^* \colon H \text{ is only linked to one or both members of }\{x_1,x_1'\}\}$. Then, $G^*$ has a component $G_v^*$ such that $\{v,w\} \subseteq V(G_v^*)$, and ${\rm C}(G^*) = \{G_v^*\} \cup \mathcal{I} \cup \mathcal{Z}$, where $\mathcal{Z}$ is $\emptyset$ or $\{(\{y^*\}, \emptyset)\}$. By Lemmas~\ref{LemIsol1} and~\ref{LemIsol2},
\begin{equation*} \iota(G,P_3) \leq 1 + \iota(G^*,P_3) = 1 + \iota(G_v^*,P_3) + \sum_{H\in \mathcal{I}} \iota(H,P_3) \leq \frac{|N[y_1]|}{4} + \iota(G_v^*,P_3) + \sum_{H\in \mathcal{I}}\frac{|V(H)|}{4}.
\end{equation*}
If $G_v^*$ is not an $\mathcal{E}$-graph, then by the induction hypothesis, $\iota(G_v^*,P_3) \leq |V(G_v^*)|/4$, so $\iota(G,P_3) \leq n/4$. Suppose that $G_v^*$ is an $\mathcal{E}$-graph. Since $d_{G^*}(v)=1$, $G_v^*$ is a $\{P_3,G_{7,1}\}$-graph. Suppose $G_v^* \simeq G_{7,1}$. Then, $d_{G_v^*}(w) = 3$ as $N_{G_v^*}(v) = \{w\}$. 
Suppose $y^* \in V(G_v^*)$. Then, $N_{G_v^*}(y^*) = \{w\}$. We have $d_{G_v^*}(v) = d_{G_v^*}(y^*) = 1$, which contradicts $G_v^* \simeq G_{7,1}$ as $G_{7,1}$ has only one vertex of degree $1$. Thus, $y^* \notin V(G_v^*)$. Let $z \in V(G_v^*)\setminus N[w]$. Then, $\{w, z\} \cup \bigcup_{H \in \mathcal{I}}D_H \in {\rm Is}(G^*, P_3)$. By Lemmas~\ref{LemIsol1} and~\ref{LemIsol2},
\begin{equation*} \iota(G,P_3) \leq 1 + \iota(G^*, P_3) \leq 3 + \sum_{H \in \mathcal{I}} \iota(H,P_3) \leq \frac{|N[y_1] \cup \{y^*\} \cup V(G_v^*)|}{4} + \sum_{H \in \mathcal{I}} \frac{|V(H)|}{4} = \frac{n}{4}.
\end{equation*}
Now suppose $G_v^* \simeq P_3$. Then, $\{w\} \cup \bigcup_{H \in \mathcal{I}}D_H \in {\rm Is}(G^*, P_3)$. If $y^* \notin V(G_v^*)$, then $\iota(G,P_3) \leq n/4$ as in the last calculation. Suppose $y^* \in V(G_v^*)$. Then, $wy^* \in E(G_v^*)$.  Recall that $N(y') = \{y_1, y^*\}$. Since $v, y^* \in N(w)$ and $\Delta(G) = 3$, $x_1^1 \notin N(w)$ for some $x_1^1 \in \{x_1, x_1'\}$. Let $x_1^2$ be the member of $\{x_1, x_1'\} \setminus \{x_1^1\}$. If $x_1^1 \notin N(y^*)$, then $G[\{v, x_1^1, y_1, y', y^*, w\}] \simeq C_6$, a contradiction. Thus, $x_1^1 \in N(y^*)$, and hence $N(y^*) = \{y', w, x_1^1\}$. If $x_1^2 \notin N(w)$, then $G[\{v, x_1^2, y_1, y', y^*, w\}] \simeq C_6$, a contradiction. Thus, $x_1^2 \in N(w)$. Since $\Delta(G) = 3$, we have $N(x_1^1) = \{v, y_1, y^*\}$, $N(x_1^2) = \{v, y_1, w\}$, and hence $\mathcal{I} = \emptyset$. Thus, we have $n = |N[y_1]| + |V(G_v^*)| = 7$, a contradiction.

Now suppose $y_1 \neq y_1'$. Let $y^*$ be the member of $V(H_1) \setminus \{y_1,y_1'\}$. By (\ref{2.2.4-1}), $d_{H_1}(y_1) \neq 2$ and $d_{H_1}(y_1') \neq 2$. Since $\Delta(H_1) = 2$, $H_1 = (\{y_1, y^*, y_1'\}, \{y_1y^*, y^*y_1'\})$. By (\ref{2.2.4-1}), $y^* \notin N[N(v)]$. 

Suppose $x_1x_1' \notin E(G)$. Since $G$ has no induced $6$-cycles, $x_1y_1' \in E(G)$ or $x_1'y_1 \in E(G)$. We may assume that $x_1'y_1 \in E(G)$. Thus, $y_1$ is adjacent to both $x_1$ and $x_1'$, and hence we obtain $\iota(G,P_3) \leq n/4$ as in the case $y_1 = y_1'$ above.

Now suppose $x_1x_1' \in E(G)$. Then, $N(x_1) = \{v, x_1', y_1\}$, $N(x_1') = \{v, x_1, y_1'\}$, and hence $\mathcal{H} = \{H_1\} \cup \mathcal{H}_w^*$. Let $G^* = G - N[x_1]$. Note that $G^*$ has a component $G_w^*$ such that $w \in V(G_w^*)$. If $wy_1' \notin E(G)$, then ${\rm C}(G^*) = \{G_w^*,G[\{y_1',y^*\}]\}$, and hence, by Lemma~\ref{LemIsol1}, Lemma~\ref{LemIsol2}, Observation~\ref{ObsE} and the induction hypothesis,
\begin{equation*} \iota(G,P_3) \leq 1 + \iota(G_w^*,P_3) \leq \frac{|N[x_1]|}{4} + \frac{|V(G_w^*)|+1}{4} < \frac{n}{4}.
\end{equation*}
Now suppose $wy_1' \in E(G)$. Since $y^* \notin N[N(v)]$, $wy^* \notin E(G)$. If $wy_1 \in E(G)$, then $N(w) = \{v, y_1, y_1'\}$, so $\mathcal{H}_w^* = \emptyset$, and hence $n = 7$, a contradiction. Thus, $wy_1 \notin E(G)$. Consequently, we have $G[\{v, x_1, y_1, y^*, y_1', w\}] \simeq C_6$, a contradiction.

Now suppose that a vertex of $H_1$ of degree $2$ in $H_1$ is adjacent to a neighbour of $v$ (in $G$). We may assume that this vertex is $y_1$. Thus, $N_{H_1}[y_1] = V(H_1)$, $N[y_1] = \{x_1\} \cup V(H_1)$, and hence $y_1 \neq y_1'$. Let $y^*$ be the member of $V(H_1) \setminus\{y_1,y_1'\}$. Let $G^* = G - N[y_1]$. Then, $G^*$ has a component $G_v^*$ such that $v, x_1', w \in V(G_v^*)$, and ${\rm C}(G^*) = \{G_v^*\} \cup \mathcal{H}_{x_1}^*$. Since $\Delta(G) = 3$ and $v, y_1 \in N(x_1)$, $|\mathcal{H}_{x_1}^*| \leq 1$. If $|\mathcal{H}_{x_1}^*| = 1$, then let $H^*$ be the member of $\mathcal{H}_{x_1}^*$. If $|\mathcal{H}_{x_1}^*| = 0$, then let $H^* = (\emptyset, \emptyset)$. By Lemma~\ref{LemIsol1}, Lemma~\ref{LemIsol2} and the induction hypothesis,
\begin{equation} \iota(G,P_3) \leq 1 + \iota(G^*,P_3) = 1 + \iota(G_v^*,P_3) + \iota(H^*,P_3) \leq \frac{|N[y_1]|}{4} + \iota(G_v^*,P_3) + \left \lfloor \frac{|V(H^*)|}{4} \right \rfloor. \label{2.2.4-2}
\end{equation}
If $G_v^*$ is not an $\mathcal{E}$-graph, then by the induction hypothesis, $\iota(G_v^*,P_3) \leq |V(G_v^*)|/4$, so $\iota(G,P_3) \leq n/4$. Suppose that $G_v^*$ is an $\mathcal{E}$-graph. By Observation~\ref{ObsE}, $\iota(G_v^*,P_3) = (|V(G_v^*)| + 1)/4$. By Observations~\ref{ObsG7} and~\ref{ObsG11}, since $d_{G_v^*}(v) = 2$ and $d_{G_v^*}(x_1') \leq |N(x_1') \setminus \{y_1'\}| \leq \Delta(G) - 1 = 2$, $G_v^*$ is a $\{P_3, C_3, C_7, C_{11}, G_{11}\}$-graph. If $G_v^*$ is a $\{C_7, C_{11}, G_{11}\}$-graph, then we obtain $\iota(G, P_3) \leq n/4$ by the arguments for Cases~2.2.2 and 2.2.3 with $y_1$ and $G_v^*$ taking the roles of $v$ and $H_1$ there, respectively. Suppose that $G_v^*$ is a $\{P_3, C_3\}$-graph. Then, $V(G_v^*) = \{v, x_1', w\}$. Since $16 \leq n = |N[y_1]| + |V(G_v^*)| + |V(H^*)| = 7 + |V(H^*)|$, $|V(H^*)| \geq 9$. Thus, $N(x_1) = \{v, y_1, z\}$ for some $z \in V(H^*)$. By the division algorithm, $|V(H^*)| = 4k + r$ for some integer $k \geq 2$ and some $r \in \{0, 1, 2, 3\}$. If $r \in \{1, 2, 3\}$, then $\floor{{|V(H^*)|}/{4}} \leq {(|V(H^*)|-1)}/{4}$, so $\iota(G, P_3) \leq n/4$ by (\ref{2.2.4-2}). Suppose $r = 0$. We have $x_1'y_1, wy_1 \notin E(G)$ as $N(y_1) = \{x_1, y_1', y^*\}$. 

Suppose $x_1'y^* \in E(G)$. Then, $N(x_1') = \{v, y_1', y^*\}$. Let $G^{\dagger} = G-N[x_1']$. Then, $V(H^{\dagger}) = V(H^*) \cup \{x_1,y_1\}$ for some $H^{\dagger} \in {\rm C}(G^{\dagger})$. Since $V(G_v^*) = \{v, x_1', w\}$, ${\rm C}(G^{\dagger}) = \{G[\{w\}],H^{\dagger}\}$. Since $|V(H^{\dagger})| = |V(H^*)| + 2 = 4k + 2$, $H^{\dagger}$ is not an $\mathcal{E}$-graph. By Lemma~\ref{LemIsol1}, Lemma~\ref{LemIsol2} and the induction hypothesis,
\begin{equation*} \iota(G,P_3) \leq 1 + \iota(G^{\dagger},P_3) \leq \frac{|N[x_1']|}{4} + \left \lfloor \frac{|V(H^{\dagger})|}{4} \right \rfloor < \frac{n}{4}.
\end{equation*}

Now suppose $x_1'y^* \notin E(G)$. Let $X' = \{x_1'\} \cup V(H_1)$. Let $G_1' = G - X'$. Then, $N_{G_1'}(v) = \{x_1, w\}$. Recall that ${\rm C}(G^*) = \{G_v^*\} \cup \mathcal{H}_{x_1}^*$ and $N_{G_v^*}[v] = \{v, x_1', w\} = V(G_v^*)$. Thus, $V(G_1') = \{v, x_1, w\} \cup V(H^*)$ and $G_1'$ is connected. Since $w \notin N(x_1)$ and no vertex of $H^*$ is adjacent to a vertex of $G_v^*$, $G_1' - v$ is not connected. By Observation~\ref{ObsE}, $G_1'$ is not a $\mathcal{G}_{15}$-graph. Since $n \geq 16$, we have $|V(G_1')| \geq 12$, so $G_1'$ is not an $\mathcal{E}$-graph. By the induction hypothesis, $\iota(G_1',P_3) \leq |V(G_1')|/4 = (n-4)/4$. 

If $d_{H_1}(y_1') = 2$, then $N[y_1'] = X'$, so $\iota(G, P_3) \leq 1 + \iota(G_1',P_3) \leq n/4$ by Lemma~\ref{LemIsol1}. 

Suppose $d_{H_1}(y_1') = 1$. Recall that $N(y_1) = \{x_1, y_1', y^*\}$, so $E(H_1) = \{y_1'y_1, y_1y^*\}$. Let $X'' = X' \setminus \{y^*\}$. Since $N(x_1) = \{v, y_1, z\}$, $x_1 \notin N(y^*)$. If $w \notin N(y^*)$, then ${\rm C}(G - X'') = \{G_1', G[\{y^*\}]\}$, and hence, since $X'' \subseteq N[y_1']$, $\iota(G,P_3) \leq 1 + \iota(G - X'', P_3) \leq 1 + \iota(G_1', P_3) \leq n/4$ by Lemmas~\ref{LemIsol1} and~\ref{LemIsol2}. Suppose $w \in N(y^*)$. Since $G[\{v, x_1', y_1', y_1, y^*, w\}] \not\simeq C_6$ (as $G$ has no induced $6$-cycles), $wy_1' \in E(G)$ or $x_1'w \in E(G)$. If $wy_1' \in E(G)$, then $N(y_1') = \{x_1', w, y_1\}$, ${\rm C}(G - N[y_1']) = \{G[\{v, x_1\} \cup V(H^*)], G[\{y^*\}]\}$, and hence, as above, 
$$\iota(G, P_3) \leq 1 + \iota(G[\{v, x_1\} \cup V(H^*)], P_3) \leq 1 + \left \lfloor \frac{4k+2}{4} \right \rfloor = 1 + k = 1 + \frac{|V(H^*)|}{4} < \frac{n}{4}.$$ 
Suppose $x_1'w \in E(G)$. Then, $N(x_1') = \{v, w, y_1'\}$. Let $G^{\dagger} = G-N[x_1']$. Then, $V(G^{\dagger}) = \{x_1,y_1,y^*\} \cup V(H^*)$ and $G^{\dagger}$ is connected. Since $n \geq 16$, $|V(G^{\dagger})| \geq 12$. By Observation~\ref{ObsG11}~(a), since $N_{G^{\dagger}}(y^*) = \{y_1\}$, $G^{\dagger}$ is not a $\mathcal{G}_{15}$-graph. By Lemma~\ref{LemIsol1} and the induction hypothesis, $\iota(G,P_3) \leq 1 + \iota(G^{\dagger},P_3) \leq 1 + \left \lfloor (4k+3)/4 \right \rfloor < n/4$. \qed
\\

We now start working towards proving Lemma~\ref{ThmP3small}.

\begin{observation}\label{ObsG7_235} Suppose that $H \in \{G_{7,2},G_{7,3},G_{7,5}\}$, $G = (V(H), E(H) \cup \{e\})$ with $e \in \binom{V(H)}{2} \setminus E(H)$, and $G$ is subcubic. If $H = G_{7,5}$ and $e \in \{\{1, 4\}, \{1,5\}\}$, then $\iota(G,P_3) \leq 1$, otherwise $G$ is a $\{G_{7,4},G_{7,6}\}$-graph.
\end{observation}

\begin{observation}\label{ObsG7_1} Suppose that $S \subseteq \binom{V(G_{7,1})}{2} \setminus E(G_{7,1})$, $1\leq |S| \leq 2$, $G = (V(G_{7,1}),E(G_{7,1}) \cup S)$, and $G$ is subcubic. If $vw \in \binom{\{2, 4, 6\}}{2}$ and $S = \{vw\}$, then $\iota(G,P_3) \leq 1$, otherwise $G$ is a $(\mathcal{G}_7\setminus\{C_7, G_{7,1}, G_{7,5}\})$-graph.
\end{observation}

\begin{lemma}\label{ThmG7}
If $G$ is a connected subcubic $n$-vertex graph with $n \leq 7$, and $G$ has no induced 6-cycles and is not a $(\mathcal{G}_3 \cup \mathcal{G}_7)$-graph, then $\iota(G,P_3) \leq n/4$.
\end{lemma}
\proof{The result follows from Lemma~\ref{LemP3G5} if $n \leq 5$. Suppose $n \geq 6$. If $\Delta(G) \leq 2$, then $\iota(G,P_3) \leq n/4$ by Lemma \ref{LemDelta2}. Suppose $\Delta(G) = 3$. Let $v \in V(G)$ with $d(v) = 3$. Let $G' = G-N[v]$. If $n = 6$, then $|V(G')| = 2$, so $\iota(G,P_3) \leq 1 < n/4$. Suppose $n = 7$. Let $m = |E(G)|$. The handshaking lemma tells us that $\sum_{z \in V(G)} d(z) = 2m$, so $3(7) \geq 2m$, and hence $m \leq 10$. The result follows if $\{z\} \in {\rm Is}(G, P_3)$ for some $z \in V(G)$. Suppose $\{v\} \notin {\rm Is}(G, P_3)$. Then, $G'$ is a $\mathcal{G}_3$-graph. Let $y_1$, $y_2$ and $y_3$ be the vertices of $G'$. Since $G$ is connected, $G'$ is linked to some $x \in N(v)$. If $G'$ is linked only to $x$, then $\{x\} \in {\rm Is}(G, P_3)$. Suppose that $G'$ is linked to some $x' \in N(v) \setminus \{x\}$. Let $w$ be the member of $N(v) \setminus \{x,x'\}$. 

Suppose that a vertex of $G'$ is adjacent to more than one neighbour of $v$. We may assume that $x, x', y_2 \in N(y_1)$. Since $\Delta(G) = 3$, $N(y_1) = \{x, x', y_2\}$. Thus, $E(G') = \{y_1y_2, y_2y_3\}$. Suppose $\{y_1\} \notin {\rm Is}(G, P_3)$. Then, $wy_3 \in E(G)$. Let $X = \{v, x', y_1, y_2, y_3, w\}$. Suppose $xx' \in E(G)$. Since $d(v) = d(x') = d(y_1) = 3 = \Delta(G)$ and $G[X] \not\simeq C_6$, $wy_2 \in E(G)$. Thus, we have $m = 10$ and $G \simeq G_{7,6}$, a contradiction. Therefore, $xx' \notin E(G)$. Suppose $x \in N(y_2)$ or $x' \in N(y_2)$. By symmetry, we may assume that $x \in N(y_2)$. Since $G[X] \not\simeq C_6$, we have $x'w \in E(G)$ or $x'y_3 \in E(G)$, and hence $m = 10$. Since $G \not\simeq G_{7,4}$, $x'w \notin E(G)$. Thus, $x'y_3 \in E(G)$, and hence $\{x'\} \in {\rm Is}(G, P_3)$. Next, suppose $w \in N(y_2)$. Since $G \not\simeq G_{7,5}$ and $d(v) = d(w) = d(y_1) = d(y_2) = 3 = \Delta(G)$ (and $xx' \notin E(G)$), we have $xy_3 \in E(G)$ or $x'y_3 \in E(G)$, so $\{z\} \in {\rm Is}(G, P_3)$ for some $z \in \{x, x'\}$. Now suppose $x, x', w \notin N(y_2)$, that is, $N(y_2) = \{y_1, y_3\}$. If $x \in N(y_3)$ or $x' \in N(y_3)$, then $\{z\} \in {\rm Is}(G, P_3)$ for some $z \in \{x, x'\}$. Suppose $x, x' \notin N(y_3)$.  From $G[X] \not\simeq C_6$ and $G[(X \setminus \{x'\}) \cup \{x\}] \not\simeq C_6$, we obtain $x'w \in E(G)$ and $xw \in E(G)$, respectively, and together with $v, y_3 \in N(w)$, this contradicts $d(w) \leq \Delta(G) = 3$. 

Now suppose that no vertex of $G'$ is adjacent to more than one neighbour of $v$. 

Suppose that a neighbour of $v$ is adjacent to two vertices of $G'$. We may assume that $y_1, y_2 \in N(x)$. Suppose $\{x\} \notin {\rm Is}(G, P_3)$. Then, $x'w \in E(G)$ and $y_3$ is adjacent to either $x'$ or $w$. We may assume that $x'y_3 \in E(G)$. Since $G \not\simeq G_{7,6}$, $G' \not\simeq C_3$. Thus, $G' \simeq P_3$, and hence $N_{G'}[y_i] = V(G')$ for some $i \in [3]$. Since $G \not\simeq G_{7,2}$, $i \notin \{1, 2\}$. We obtain $i = 3$, which yields $G \simeq G_{7,5}$, a contradiction.  

Now suppose that no neighbour of $v$ is adjacent to more than one vertex of $G'$. Since $G'$ is linked to $x$, we may assume that $xy_1 \in E(G)$. Suppose $\{x\} \notin {\rm Is}(G, P_3)$. Then, $uy_i \in E(G)$ for some $u \in N(v) \setminus \{x\}$ and $i \in \{2, 3\}$. We may assume that $uy_i = x'y_2$. 

Suppose $y_1y_2 \in E(G)$. Then, $d_{G'}(y_1) = 2$ or $d_{G'}(y_2) = 2$. We may assume that $d_{G'}(y_2) = 2$. Thus, $y_2y_3 \in E(G)$. Recall that $m \leq 10$. Suppose $xw \in E(G)$. Then, $G$ contains a copy $G^*$ of $G_{7,1}$. Since $G$ is not a $\mathcal{G}_7$-graph, we have $|E(G^*)| + 1 \leq m \leq 10 = |E(G^*)| + 2$, so $\iota(G, P_3) \leq 1$ by Observation~\ref{ObsG7_1}. Similarly, $\iota(G, P_3) \leq 1$ if $y_1y_3 \in E(G)$. If $xw, y_1y_3 \notin E(G)$, then $\{x'\} \in {\rm Is}(G, P_3)$.

Now suppose $y_1y_2\notin E(G)$. Then, $G' \simeq P_3$ and $E(G') = \{y_1y_3,y_2y_3\}$. Since $G[\{v, x, y_1, y_3, y_2, x'\}] \not\simeq C_6$, $xx' \in E(G)$. Suppose $\{x\} \notin {\rm Is}(G, P_3)$. Then, $wy_3 \in E(G)$. Thus, $G$ contains a copy of $G_{7,3}$, and hence we have $10 = |E(G_{7,3})| + 1 \leq m \leq 10$. By Observation~\ref{ObsG7_235}, $G$ is a $\mathcal{G}_7$-graph, a contradiction. \qed
}

\begin{lemma}\label{ThmG8-11}
If $G$ is a connected subcubic $n$-vertex graph with $8 \leq n \leq 11$, and $G$ has no induced $6$-cycles and is not a $\mathcal{G}_{11}$-graph, then $\iota(G,P_3) \leq 2$.
\end{lemma}
\proof{If $\Delta(G) \leq 2$, then $\iota(G,P_3) \leq 2$ by Lemma~\ref{LemDelta2}. Suppose $\Delta(G) = 3$. Let $v \in V(G)$ with $d(v) = 3$, let $G' = G - N[v]$, and let $n' = |V(G')|$. By Lemma~\ref{LemIsol1}, $\iota(G,P_3) \leq 1 + \iota(G',P_3)$. Thus, the result follows if $\iota(G',P_3) \leq 1$. Define $\mathcal{H}$ and $\mathcal{H}'$ as in the proof of Theorem~\ref{ThmP3}. We assume that, as in Case~2.2 of the proof of Theorem~\ref{ThmP3}, $|\mathcal{H}'| = 1$ and the member $H_1$ of $\mathcal{H}'$ is linked to some $x, x' \in N(v)$ with $x \neq x'$, because otherwise the result follows by the argument in the proof of Theorem~\ref{ThmP3} (using Lemma~\ref{ThmG7} instead of the induction hypothesis). Let $w$ be the member of $N(v) \setminus \{x, x'\}$, and let $y_1 = y_{x, H_1}$ and $y_1' = y_{x', H_1}$. By Lemma~\ref{ThmG7}, since $n' = n - 4 \leq 7$, $H_1$ is a $(\mathcal{G}_3 \cup \mathcal{G}_7)$-graph, no member of $\mathcal{H} \setminus \{H_1\}$ is an $\mathcal{E}$-graph, and if $\mathcal{H} \setminus \{H_1\} \neq \emptyset$, then $H_1$ is a $\mathcal{G}_3$-graph. Since each connected $3$-vertex graph is a $\mathcal{G}_3$-graph, no member of $\mathcal{H} \setminus \{H_1\}$ is of order~$3$.\medskip

\noindent\textbf{Case 1:} \emph{$\mathcal{H} \setminus \{H_1\}$ has a member $H_2$.} Then, $|V(H_1)| = 3$ and $|V(H_2)| \in \{1, 2, 4\}$. If $|V(H_2)| \in \{1, 2\}$, then $|V(H)| \leq 2$ for each $H \in \mathcal{H} \setminus \{H_1\}$, so $\{v, y_1\} \in {\rm Is}(G, P_3)$. Suppose $|V(H_2)| = 4$. Then, $\mathcal{H} = \{H_1, H_2\}$. Suppose that $H_2$ is linked only to some $u \in N(v)$. Since $H_1$ is linked to $x$ and $x'$, ${\rm C}_3(G-N[y_{u,H_2}])$ has only one member $F$, where $V(F) = V(H_1) \cup (N[v] \setminus \{u\})$. Thus, $\iota(G,P_3) \leq 2$ by Lemmas~\ref{LemIsol1} and~\ref{ThmG7}. Now suppose that $H_2$ is linked to at least two neighbours of $v$. 

Suppose that some $y \in V(H_1)$ is adjacent to two neighbours of $v$. We may assume that $y = y_1 = y_1'$, meaning that $N(y_1) = \{x, x', y_2\}$ for some $y_2 \in V(H_1) \setminus \{y_1\}$. Thus, $E(H_1) = \{y_1y_2, y_2y_3\}$, where $y_3$ is the member of $V(H_1) \setminus \{y_1, y_2\}$. If $H_2$ is linked to $w$, then $\{y_1, y_{w, H_2}\} \in {\rm Is}(G, P_3)$. Suppose that $H_2$ is linked only to $x$ and $x'$. Then, $N(x) = \{v, y_1, y_{x, H_2}\}$ and $N(x') = \{v, y_1, y_{x', H_2}\}$. If $wy_3 \notin E(G)$, then ${\rm C}_3(G - N[y_1]) = \{H_2\}$, so $\{y_1, y_{x, H_2}\} \in {\rm Is}(G, P_3)$. Suppose $wy_3 \in E(G)$. Since $G[\{v, x, y_1, y_2, y_3, w\}] \not\simeq C_6$, $wy_2 \in E(G)$. Let $G^* = G - N[w]$. Then, $G^* = G[\{x,x',y_1\} \cup V(H_2)]$, and hence $G^*$ is a connected $7$-vertex graph. Suppose $\iota(G^*,P_3) \geq 2$. By Lemma~\ref{ThmG7}, $G^*$ is a $\mathcal{G}_7$-graph. Since $d_{G^*}(x) = d_{G^*}(y_1) = d_{G^*}(x') = 2$, $G^* \simeq C_7$ by Observation~\ref{ObsG7}. We obtain that $G \simeq G_{11}$, a contradiction. Thus, $\iota(G^*,P_3) \leq 1$. By Lemma~\ref{LemIsol1}, $\iota(G,P_3) \leq 1 + \iota(G^*,P_3) \leq 2$.  

Now suppose that no vertex of $H_1$ is adjacent to more than one neighbour of $v$. Then, $y_1 \neq y_1'$. Let $y^*$ be the member of  $V(H_1) \setminus \{y_1,y_1'\}$.

Suppose that some $u \in N(v)$ is adjacent to two vertices of $H_1$. We may assume that $u = x$. Thus, $N(x) = \{v, y_1, y^*\}$, and hence $H_2$ is linked to $x'$ and $w$ (but not to $x$). Thus, $N(x') = \{v, y_1', y_{x',H_2}\}$. Since each vertex of $H_1$ is adjacent to $x$ or $x'$, no vertex of $H_1$ is adjacent to $w$. If two vertices of $H_2$ are adjacent to $w$, then $\{w, y_1'\} \in {\rm Is}(G, P_3)$. Suppose $N(w) \cap V(H_2) = \{y_{w,H_2}\}$. Since $H_2$ is connected, $z_1 \in N(y_{x',H_2})$ for some $z_1 \in V(H_2)$. Let $I = G[(V(H_2) \cup \{w\})\setminus N[y_{x',H_2}]]$. If $d(y_{x',H_2}) = 3$ or $I$ is not connected, then $\{x,y_{x',H_2}\} \in {\rm Is}(G, P_3)$. Now suppose that $N(y_{x',H_2}) = \{x', z_1\}$ and $I$ is connected. Then, since $N(w) \cap V(H_2) = \{y_{w,H_2}\}$, we have $y_{w,H_2} \notin N[y_{x',H_2}]$ and $E(I) = \{wy_{w,H_2}, y_{w,H_2}z_2\}$, where $\{z_2\} = V(H_2) \setminus \{y_{x',H_2}, z_1, y_{w,H_2}\}$. If $z_1 \in N(y_{w,H_2})$, then $G[\{v, x', y_{x',H_2}, z_1, y_{w,H_2}, w\}] \simeq C_6$, a contradiction. Thus, $z_1 \notin N(y_{w,H_2})$. Therefore, since $H_2$ is connected, $E(H_2) = \{y_{x',H_2}z_1, z_1z_2, z_2y_{w,H_2}\}$. Thus, $\{x, z_1\} \in {\rm Is}(G, P_3)$. 

Now suppose that no neighbour of $v$ is adjacent to more than one vertex of $H_1$. Thus, $N(x) \cap V(H_1) = \{y_1\}$ and $N(x') \cap V(H_1) = \{y_1'\}$. Suppose $y_1y_1' \notin E(G)$. Then, $E(H_1) = \{y_1y^*, y^*y_1'\}$. Since $\Delta(G) = 3$ and $H_2$ is linked to at least one of $x$ and $x'$, we have $xx' \notin E(G)$, and hence $G[\{v, x, y_1, y^*, y_1', x'\}] \simeq C_6$, a contradiction. Thus, $y_1y_1' \in E(G)$. Since $H_1$ is a $\mathcal{G}_3$-graph, $d_{H_1}(y_1) = 2$ or $d_{H_1}(y_1') = 2$. We may assume that $d_{H_1}(y_1) = 2$. Thus, $N[y_1] = V(H_1) \cup \{x\}$. Let $G^* = G - N[y_1]$. Since $H_2$ is linked to at least two neighbours of $v$, $G^*$ is connected. By Lemma~\ref{LemIsol1}, $\iota(G,P_3) \leq 1 + \iota(G^*,P_3)$. Thus, $\iota(G,P_3) \leq 2$ if $\iota(G^*,P_3) \leq 1$. Suppose $\iota(G^*,P_3) \geq 2$. Since $d_{G^*}(v) = 2$ and $d_{G^*}(x') \leq 2$, $G^* \simeq C_7$ by Observation \ref{ObsG7}. Thus, $H_2 \simeq P_4$, $H_2$ is linked to $x'$ and $w$, $y_{x',H_2} \neq y_{w,H_2}$ and $E(H_2) = \{y_{x',H_2}z_1, z_1z_2, z_2y_{w,H_2}\}$, where $\{z_1, z_2\} = V(H_2) \setminus \{y_{x',H_2}, y_{w,H_2}\}$. We have $N(x') = \{v, y_1', y_{x',H_2}\}$. 

Suppose $y_1'y^* \notin E(G)$. Recall that no vertex of $H_1$ is adjacent to more than one neighbour of $v$. Thus, $N(y_1') \cap N(v) = \{x'\}$. If $y^* \in N(w)$, then $N(w) = \{v, y_{w,H_2}, y^*\}$ and $G[\{v, x', y_1', y_1, y^*, w\}] \simeq C_6$, a contradiction. Thus, $y_1', y^* \notin N(w)$, and hence $\{x, z_1\} \in {\rm Is}(G, P_3)$.  

Now suppose $y_1'y^* \in E(G)$. If $N(x) \cap \{z_1, y_{x',H_2}\} = \emptyset$, then $\{y_1',y_{w,H_2}\} \in {\rm Is}(G, P_3)$. Suppose $N(x) \cap \{z_1, y_{x',H_2}\} \neq \emptyset$. If $z_1 \in N(x)$, then $N(x) = \{v, y_1, z_1\}$, and hence $G[\{v, x, z_1, z_2, y_{w,H_2}, w\}] \simeq C_6$ (as $G^* \simeq C_7$), a contradiction. Thus, $y_{x',H_2} \in N(x)$, and hence $N(x) = \{v, y_1, y_{x',H_2}\}$ and $N(y_{x',H_2}) = \{x, x', z_1\}$. If $wy^* \notin E(G)$, then $G \simeq G_{11}$, a contradiction. Thus, $wy^* \in E(G)$, and hence $\{y^*, y_{x',H_2}\} \in {\rm Is}(G, P_3)$.\medskip

\noindent\textbf{Case 2:} \emph{$\mathcal{H} \setminus \{H_1\} = \emptyset$.} We have $|V(H_1)| = n' = n - 4 \geq 4$, so $H_1$ is a $\mathcal{G}_7$-graph. 
We may assume that $H_1 \in \mathcal{G}_7$. Thus, $V(H_1) = [7]$. We keep in mind that since $\Delta(G) = 3$, $N(i) = N_{H_1}(i)$ for each $i \in V(H_1)$ with $d_{H_1}(i) = 3$.

Suppose that some vertex of $H_1$ is adjacent to two neighbours of $v$. Thus, we may assume that $y_1 = y_1'$, meaning that $N(y_1) = \{x, x', y\}$ for some $y \in V(H_1) \setminus \{y_1\}$. Since $d_{H_1}(y_1) = 1$, $H_1 = G_{7,1}$. Thus, $y_1$ is the vertex $1$ of $G_{7,1}$. For each $i \in [7] \setminus \{1\}$, let $y_i$ be the vertex $i$ of $H_1$. 
Then, $y = y_7$. If $wy_2 \notin E(G)$, then $\{y_1, y_5\} \in {\rm Is}(G, P_3)$. If $wy_6 \notin E(G)$, then $\{y_1,y_3\} \in {\rm Is}(G, P_3)$. If $wy_2, wy_6 \in E(G)$, then $N(w) = \{v, y_2, y_6\}$, $N(y_2) = \{w, y_3, y_7\}$ and $G[\{v, x, y_1, y_7, y_2, w\}] \simeq C_6$, a contradiction.  

Now suppose $|N(z) \cap N(v)| \leq 1$ for each $z \in V(H_1)$. Then, $y_1 \neq y_1'$.\medskip

\noindent\textbf{Case 2.1:} \emph{$H_1 \neq C_7$.} By Observation~\ref{ObsG7}~(a), there exists some $y^* \in V(H_1) \setminus N[y_1]$ such that $d_{H_1}(y^*) = 3 = d(y^*)$ and $H_1 - N[y^*]$ is connected. Let $H^* = H_1 - N[y^*]$. We have that $y_1 \in V(H^*)$, $H^*$ is a $\{P_3, C_3\}$-graph, and $G - N[y^*]$ is a connected $7$-vertex graph (as $xy_1 \in E(G)$). Let $z_1 \in N_{H^*}(y_1)$, and let $z_2$ be the member of $V(H^*) \setminus \{y_1, z_1\}$. Since $x \in N(y_1)$ and $\Delta(G) = 3$, $d_{H_1}(y_1) \leq 2$. By Observation~\ref{ObsG7}~(b), $d_{H_1}(z) = 3 = d(z)$ for each $z \in N_{H_1}(y_1)$. Thus, $N(z_1) = N_{H_1}(z_1)$, meaning that $N(z_1) \cap N(v) = \emptyset$. If $z_2 \in N_{H_1}(y_1)$, then $N(z_2) = N_{H_1}(z_2)$, so $\{x,y^*\} \in {\rm Is}(G, P_3)$. Suppose $z_2 \notin N_{H_1}(y_1)$. Then, $E(H^*) = \{y_1z_1, z_1z_2\}$. If $u \notin N(z_2)$ for each $u \in N(v) \setminus \{x\}$, then $\{x, y^*\} \in {\rm Is}(G, P_3)$. Suppose $u \in N(z_2)$ for some $u \in N(v) \setminus \{x\}$. Recall that $|N(z) \cap N(v)| \leq 1$ for each $z \in V(H_1)$. Thus, $N(y_1) \cap N(v) = \{x\}$, $N(z_2) \cap N(v) = \{u\}$, and since $G[\{v, x, y_1, z_1, z_2, u\}] \not\simeq C_6$, $xu \in E(G)$. Consequently, $\{x, y^*\} \in {\rm Is}(G, P_3)$.\medskip

\noindent\textbf{Case 2.2:} \emph{$H_1 = C_7$.} We may assume that $y_1 = 1$. Let $y_i = i$ for each $i \in [7] \setminus \{1\}$. Then, $N(y_1) = \{x, y_2, y_7\}$. Let $I = G - N[y_1]$. Then, $|V(I)| = 7$.

Suppose that $I$ is connected. If $I$ is not a $\mathcal{G}_7$-graph, then $\iota(G,P_3) \leq 1 + \iota(I,P_3) \leq 2$ by Lemmas~\ref{LemIsol1} and~\ref{ThmG7}. Suppose that $I$ is a $\mathcal{G}_7$-graph. If $I \not\simeq C_7$, then the result follows as in Case~2.1 (with $y_1$ here taking the role of $v$ there). Suppose $I \simeq C_7$. Then, $E(I) = \{vu_1, u_1y_3, y_3y_4, y_4y_5, y_5y_6, y_6u_2, u_2v\}$ with $\{u_1, u_2\} = \{x', w\}$. We have $d(v) = d(y_1) = d(y_3) = d(y_6) = 3 = \Delta(G)$ and $v, y_1 \in N(x)$. Since $G[\{v, x, y_1, y_2, y_3, u_1\}] \not\simeq C_6$ and $G[\{v, x, y_1, y_7, y_6, u_2\}] \not\simeq C_6$, we have $xy_2, u_2y_7 \in E(G)$ or $xu_1, u_2y_7 \in E(G)$ or $u_1y_2, xy_7 \in E(G)$ or $u_1y_2, xu_2 \in E(G)$ or $u_1y_2, u_2y_7 \in E(G)$. 
If $xy_2, u_2y_7 \in E(G)$, then $\{u_2, y_3\} \in {\rm Is}(G, P_3)$. If $xu_1, u_2y_7 \in E(G)$, then $\{x, y_5\} \in {\rm Is}(G, P_3)$. If $u_1y_2, xy_7 \in E(G)$, then $\{u_1, y_6\} \in {\rm Is}(G, P_3)$. If $u_1y_2, xu_2 \in E(G)$, then $\{x, y_4\} \in {\rm Is}(G, P_3)$. Suppose $u_1y_2, u_2y_7 \in E(G)$. Then, $N(u_1) = \{v, y_2, y_3\}$, $N(u_2) = \{v, y_6, y_7\}$, $N(y_2) = \{u_1, y_1, y_3\}$ and $N(y_7) = \{u_2, y_1, y_6\}$. We have $G[\{v, u_1, y_2, y_1, y_7, u_2\}] \simeq C_6$, a contradiction.

Now suppose that $I$ is not connected. Thus, $y_1' \in \{y_2, y_7\}$ and $N(x') \cap V(H_1), N(w) \cap V(H_1) \subseteq \{y_2, y_7\}$ (recall that $N(y_1) = \{x, y_2, y_7\}$). We may assume that $y_1' = y_2$. Since $\Delta(G) = 3$, $N(y_2) = \{x',y_1,y_3\}$, so $N(w) \cap V(H_1) \subseteq \{y_7\}$. If $wy_7 \in E(G)$, then $G - N[y_2]$ is connected, so the result follows as in the case where $I$ is connected (with $G-N[y_2]$ here taking the role of $I$ there). Suppose $wy_7 \notin E(G)$. Then, $N(w) \cap V(H_1) = \emptyset$. If $x'y_7 \in E(G)$, then since $G \not\simeq G_{11}$, we have $xw \notin E(G)$, so at least one of $\{x', y_4\}$ and $\{x', y_5\}$ is a $P_3$-isolating set of $G$ (as $|N(x) \cap \{y_3, y_4, y_5, y_6\}| \leq \Delta(G) - |N(x) \cap \{v, y_1\}| = 1$). Suppose $x'y_7 \notin E(G)$. Then, $N(x') \cap V(H_1) = \{y_2\}$. If $xy_4 \in E(G)$, then $G - N[y_4]$ is connected, so the result follows as in the case where $I$ is connected. If $xy_3 \in E(G)$, then since $G \not\simeq G_{11}$, we have $x'w \notin E(G)$, so $\{x, y_5\} \in {\rm Is}(G, P_3)$. Suppose $y_3, y_4 \notin N(x)$. If $xw \notin E(G)$, then $\{x', y_6\} \in {\rm Is}(G, P_3)$. If $xw \in E(G)$, then $\{x, y_4\} \in {\rm Is}(G, P_3)$. \qed}

\begin{lemma}\label{ThmG15}
If $G$ is a connected subcubic $n$-vertex graph with $12 \leq n \leq 15$, and $G$ has no induced $6$-cycles and is not a copy of $G_{15}$, then $\iota(G,P_3) \leq 3$.
\end{lemma}

\proof{If $\Delta(G) \leq 2$, then $\iota(G,P_3) \leq 3$ by Lemma~\ref{LemDelta2}. Suppose $\Delta(G) = 3$. Let $v \in V(G)$ with $d(v) = 3$, let $G' = G - N[v]$, and let $n' = |V(G')|$. If $\iota(G',P_3) \leq 2$, then $\iota(G,P_3) \leq 3$ by Lemma~\ref{LemIsol1}. Define $\mathcal{H}$ and $\mathcal{H}'$ as in the proof of Theorem~\ref{ThmP3}. We assume that, as in Case~2.2 of the proof of Theorem~\ref{ThmP3}, $|\mathcal{H}'| = 1$ and the member $H_1$ of $\mathcal{H}'$ is linked to some $x, x' \in N(v)$ with $x \neq x'$, because otherwise the result follows by the argument in the proof of Theorem~\ref{ThmP3} (using Lemmas~\ref{ThmG7} and~\ref{ThmG8-11} instead of the induction hypothesis). Let $w$ be the member of $N(v) \setminus \{x, x'\}$, and let $y = y_{x, H_1}$ and $y' = y_{x', H_1}$. By Lemmas~\ref{ThmG7} and~\ref{ThmG8-11}, since $n' = n - 4 \leq 11$, $H_1$ is a $(\mathcal{G}_3 \cup \mathcal{G}_7 \cup \mathcal{G}_{11})$-graph, and no member of $\mathcal{H} \setminus \mathcal{H}'$ is an $\mathcal{E}$-graph. Since each connected $3$-vertex graph is a $\mathcal{G}_3$-graph, no member of $\mathcal{H}\setminus\mathcal{H}'$ is of order $3$. 

Suppose $12 \leq n \leq 14$. Then, $n' \leq 10$, so $H_1$ is a $(\mathcal{G}_3 \cup \mathcal{G}_7)$-graph. If $|V(H_1)| = 7$, then $|V(H)| \leq 2$ for each $H \in \mathcal{H} \setminus \mathcal{H}'$, so $\iota(G', P_3) = \iota(H_1, P_3) = 2$ by Lemma~\ref{LemIsol2}. Suppose $|V(H_1)| = 3$. Then, $n' - |V(H_1)| \leq 7$, so at most one member of $\mathcal{H} \setminus \mathcal{H}'$ is of order at least $4$. If such a member exists, then its $P_3$-isolation number is $1$ by Lemma~\ref{ThmG7}. Since no member of $\mathcal{H} \setminus \mathcal{H}'$ is of order $3$, it follows that $\iota(G', P_3) \leq \iota(H_1, P_3) + 1 = 2$ by Lemma~\ref{LemIsol2}. 

Now suppose $n = 15$. We keep in mind that since $\Delta(G) = 3$, $N(z) = N_{H_1}(z)$ for each $z \in V(H_1)$ with $d_{H_1}(z) = 3$, and $d_{H_1}(z) \leq |N(z) \setminus N(v)| \leq 2$ for each $z \in V(H_1)$ with $N(z) \cap N(v) \neq \emptyset$. Thus, $d_{H_1}(y) \leq 2$ and $d_{H_1}(y') \leq 2$.\medskip

\noindent\textbf{Case~1:} \emph{$H_1$ is a $\mathcal{G}_{11}$-graph.} We may assume that $H_1 \in \mathcal{G}_{11}$. Thus, $V(H_1) = [11]$. Let $y_i = i$ for each $i \in [11]$.\medskip 

\noindent\textbf{Case~1.1:} \emph{$H_1 = G_{11}$.} 

Suppose $N(y_1) \cap N(v) \neq \emptyset$. We may assume that $y_1 = y$. Thus, $N(y_1) = \{x, y_2, y_{11}\}$, $d_{H_1}(y_i) = d(y_i) = 3$ for each $i \in \{2, 3, 4, 9, 10, 11\}$, $y' \in \{y_5, y_6, y_7, y_8\}$, and hence $G - N[y_3]$ is a connected $11$-vertex graph with $d_{G-N[y_3]}(y_{10}) = 1$. By Lemmas~\ref{LemIsol1} and~\ref{ThmG8-11}, and Observation~\ref{ObsG11}~(a), $\iota(G, P_3) \leq 1 + \iota(G-N[y_3], P_3) \leq 3$.

Now suppose $N(y_1) \cap N(v) = \emptyset$. Thus, $y, y' \in \{y_5, y_6, y_7, y_8\}$. Let $G^* = G - N[y_2]$. Since $H_1 - N_{H_1}[y_2] \simeq C_7$, $G^*$ is a connected $11$-vertex graph. If $\iota(G^*,P_3) \leq 2$, then $\iota(G,P_3) \leq 3$ by Lemma~\ref{LemIsol1}. Suppose $\iota(G^*,P_3) > 2$. By Lemma~\ref{ThmG8-11}, $G^*$ is a $\mathcal{G}_{11}$-graph. Since $d(y) = 3 = d_{G^*}(y)$, we have $G^* \not\simeq C_{11}$, so $G^* \simeq G_{11}$. For each $i \in \{4, 9, 10\}$, $d_{G^*}(y_i) = d(y_i) - |N(y_i) \cap N(y_2)| = 2$. Thus, since $G^* \simeq G_{11}$, we have $\{5,8\} = \{j, h\}$ with $d_{G^*}(y_j) = 2$ and $d_{G^*}(y_h) = 3$. 
Suppose $j = 5$. Then, for each $i \in \{6, 7, 8\}$, $d_{G^*}(y_i) = 3$, so $y_iu \in E(G^*)$ for some $u \in N(v)$. We may assume that $y = y_8$. Since $G^* \simeq G_{11}$, we obtain $N(x) = \{v, y_6, y_8\}$, $N(x') = \{v, w, y_7\}$ and $N(w) = \{v, x'\}$. This yields $G \simeq G_{15}$, a contradiction. Thus, $j = 8$. 
Since $G^* \simeq G_{11}$, for each $i \in \{5, 6, 7\}$, $d_{G^*}(y_i) = 3$, so $y_iu \in E(G^*)$ for some $u \in N(v)$. We may assume that $y = y_5$. We obtain $N(x) = \{v, y_5, y_7\}$, $N(x') = \{v, w, y_6\}$ and $N(w) = \{v, x'\}$. This yields $G \simeq G_{15}$, a contradiction.\medskip

\noindent\textbf{Case~1.2:} \emph{$H_1 = C_{11}$.} We may assume that $y = y_1$. Thus, $N(y_1) = \{x, y_2, y_{11}\}$. 

Suppose that $w$ is not linked to $H_1$. Then, ${\rm C}_3(G-N[\{y, y'\}]) = \{H^*\}$ with $H^* \simeq P_i$ for some $i \in [7] \setminus [2]$. Therefore, by Lemmas~\ref{LemIsol1} and~\ref{ThmG7}, $\iota(G,P_3) \leq 2 + \iota(H^*,P_3) \leq 3$. 

Now suppose that $w$ is linked to $H_1$. Let $y'' = y_{w,H_1}$. Let $Y = \{y_1, y', y''\}$. For each $z \in Y$, $|N(z) \cap N(v)| = 1$ as $\Delta(G) = 3$ and $d_{H_1}(z) = 2$. Thus, since $xy_1, x'y', wy'' \in E(G)$, $y_1$, $y'$ and $y''$ are distinct. Since $H_1 = C_{11}$, $d_{H_1}(z, z') \in [5]$ for any $z, z' \in V(H_1)$ with $z \neq z'$.\medskip

\noindent\textbf{Case~1.2.1:} \emph{$d_{H_1}(z, z') \in \{3,4,5\}$ for every $z, z' \in Y$ with $z \neq z'$.}  Then, $N_{H_1}(z) \cap N_{H_1}(z') = \emptyset$ for every $z, z' \in Y$ with $z \neq z'$. Since $|N_{H_1}[y_i]| = 3$ for each $i \in [11]$, $|V(H_1 - N_{H_1}[Y])| = 2$. Since $N(v) \subseteq N[Y]$, $Y \in {\rm Is}(G, P_3)$.\medskip

\noindent\textbf{Case~1.2.2:} \emph{$d_{H_1}(z, z') = 2$ for some $z, z' \in Y$.} We may assume that $z = y_1$ and $z' = y' = y_3$. Thus, $N(y_3) = \{x', y_2, y_4\}$. Since $G[\{v, x, y_1, y_2, y_3, x'\}] \not\simeq C_6$, we have $xy_2\in E(G)$ or $x'y_2\in E(G)$ or $xx' \in E(G)$.

Suppose $xy_2 \in E(G)$ or $x'y_2 \in E(G)$. By symmetry, we may assume that $xy_2 \in E(G)$. Thus, $N(x) = \{v, y_1, y_2\}$, $N(y_2) = \{x, y_1, y_3\}$ and $y'' \notin \{y_1, y_2\}$ (as $w \in N(y'')$). Let $G^* = G - N[x]$. Thus, $|V(G^*)| = 11$. If $\iota(G^*,P_3) \leq 2$, then $\iota(G,P_3) \leq 3$ by Lemma~\ref{LemIsol1}. Suppose $\iota(G^*,P_3) > 2$. Since $x'y_3, wy'' \in E(G^*)$, $G^*$ is connected. By Lemma~\ref{ThmG8-11}, $G^*$ is a $\mathcal{G}_{11}$-graph. We have $d_{G^*}(x') \leq 2$, $d_{G^*}(w) \leq 2$ and $N_{G^*}(y_3) = \{x', y_4\}$. Since $|(N_{G^*}(x') \cup N_{G^*}(w)) \cap (V(G^*) \setminus \{x', w, y_3\})| \leq |N_{G^*}(x') \setminus \{y_3\}| + d_{G^*}(w) \leq 1 + 2 = 3$, $G^*$ has at most $3$ vertices of degree $3$ in $G^*$ (all of which are members of $V(H_1)\setminus \{y_1, y_2, y_3\}$). Thus, $G^* \not\simeq G_{11}$, and hence $G^* \simeq C_{11}$. Consequently, $E(G^*) = \{x'y_3, y_3y_4, \dots, y_{10}y_{11}, y_{11}w, wx'\}$, $N(x') = \{v, w, y_3\}$, $G-N[x']$ is a connected $11$-vertex graph, and $d_{G-N[x']}(y_4) = 1$. By Observation~\ref{ObsG11} (a) and Lemma~\ref{ThmG8-11}, $\iota(G-N[x'], P_3) \leq 2$. By Lemma~\ref{LemIsol1}, $\iota(G,P_3) \leq 1 + \iota(G-N[x'], P_3) \leq 3$. 

Now suppose $xx' \in E(G)$. Then, $N(x) = \{v, x', y_1\}$ and $N(x') = \{v, x, y_3\}$. Let $G^* = G - N[x]$. Thus, $|V(G^*)| = 11$. If $\iota(G^*,P_3) \leq 2$, then $\iota(G,P_3) \leq 3$ by Lemma~\ref{LemIsol1}. Suppose $\iota(G^*,P_3) > 2$. Since $wy'' \in E(G^*)$, $G^*$ is connected. By Lemma~\ref{ThmG8-11}, $G^*$ is a $\mathcal{G}_{11}$-graph. We have $d_{G^*}(w) \leq 2$ and $d_{G^*}(y_3) = 2$. Since $|N_{G^*}(w) \cap (V(G^*) \setminus \{y_1, y_3\})| \leq 2$, $G^*$ has at most $2$ vertices of degree $3$ in $G^*$ (each of which is a member of $V(H_1) \setminus \{y_1, y_3\})$.  Thus, $G^* \not\simeq G_{11}$, and hence $G^* \simeq C_{11}$. Consequently, $E(G^*) = \{wy_2, y_2y_3, \dots, y_{10}y_{11}, y_{11}w\}$, $G-N[x']$ is a connected $11$-vertex graph, and $d_{G-N[x']}(y_4) = 1$. By Observation~\ref{ObsG11} (a) and Lemma~\ref{ThmG8-11}, $\iota(G-N[x'], P_3) \leq 2$. By Lemma~\ref{LemIsol1}, $\iota(G,P_3) \leq 1 + \iota(G-N[x'], P_3) \leq 3$.\medskip

\noindent\textbf{Case~1.2.3:} \emph{$d_{H_1}(u, u') \neq 2$ for every $u, u' \in Y$ with $u \neq u'$, and $d_{H_1}(z, z') = 1$ for some $z, z' \in Y$.} We may assume that $z = y_1$ and $z' = y' = y_2$. We have $y'' \neq y_1$, $y'' \neq y_2$, $d_{H_1}(y'', y_1) \neq 2$ and $d_{H_1}(y'', y_2) \neq 2$, so $y'' \notin \{y_1, y_2, y_3, y_4, y_{10}, y_{11}\}$. If $y'' = y_7$, then $Y \in {\rm Is}(G, P_3)$. Suppose $y'' \neq y_7$. Then, $y'' \in \{y_5, y_6, y_8, y_9\}$. Let $G^* = G-N[y_1]$. Note that $G^*$ is a connected $11$-vertex graph. If $d_{G^*}(y_{10}) = 1$, then $\iota(G,P_3) \leq 1 + \iota(G^*,P_3) \leq 3$ by Lemma~\ref{LemIsol1}, Observation~\ref{ObsG11} (a) and Lemma~\ref{ThmG8-11}. Suppose $d_{G^*}(y_{10}) > 1$. Then, $x'y_{10} \in E(G)$ or $wy_{10} \in E(G)$, and hence, since $d_{H_1}(y_1, y_{10}) = 2$, the result is obtained from Case~1.2.2 as follows. If $x'y_{10} \in E(G)$, then we can take $y'$ to be $y_{10}$ instead of $y_2$. If $wy_{10} \in E(G)$, then we can take $y''$ to be $y_{10}$ instead.\medskip 

\noindent\textbf{Case~2:} \emph{$H_1$ is a $\mathcal{G}_7$-graph.} We may assume that $H_1 \in \mathcal{G}_7$. Thus, $V(H_1) = [7]$. Let $y_i = i$ for each $i \in [7]$. Suppose $H_1 \neq C_7$. By Observation~\ref{ObsG7} (a), there exists some $y^* \in V(H_1)\setminus N[y]$ such that $d_{H_1}(y^*) = 3$ and $H_1 - N_{H_1}[y^*]$ is connected. Thus, $G - N[y^*]$ is a connected $11$-vertex graph. If $G - N[y^*]$ is a $\mathcal{G}_{11}$-graph, then the result follows as in Case~1 (with $y^*$ here taking the role of $v$ there). If $G - N[y^*]$ is not a $\mathcal{G}_{11}$-graph, then $\iota(G,P_3) \leq 1 + \iota(G - N[y^*], P_3) \leq 3$ by Lemmas~\ref{LemIsol1} and~\ref{ThmG8-11}. Now suppose $H_1 = C_7$. We may assume that $y = y_1$. Thus, $N(y_1) = \{x, y_2, y_7\}$. Since $x'y' \in E(G)$, $y' \neq y_1$.

If no member of $\mathcal{H} \setminus \mathcal{H}'$ contains a $P_3$-copy, then $\iota(G, P_3) \leq 1 + \iota(G - N[v], P_3) = 1 + \iota(H_1, P_3) = 3$ by Lemma~\ref{LemIsol1}, Lemma~\ref{LemIsol2} and Observation~\ref{ObsE} (a). Suppose that a member of $\mathcal{H} \setminus \mathcal{H}'$ contains a $P_3$-copy. Since $|V(G') \setminus V(H_1)| = 4$ and no member of $\mathcal{H} \setminus \mathcal{H}'$ is of order $3$, $\mathcal{H} \setminus \mathcal{H}' = \{H^*\}$ for some connected $4$-vertex graph $H^*$. If $d_{H^*}(z) = 3$ for some $z \in V(H^*)$, then $N[z] = N_{H^*}[z] = V(H^*)$, $G - N[z]$ is the $11$-vertex graph $G[N[v] \cup V(H_1)]$ and is connected, and hence the result follows as in Case~1 (or as in the proof of Theorem~\ref{ThmP3} if $|N[V(G - N[z])] \cap N(z)| = 1$) or by Lemmas~\ref{LemIsol1} and~\ref{ThmG8-11}. Suppose $d_{H^*}(z) \leq 2$ for each $z \in V(H^*)$. Then, $H^* \simeq P_4$ or $H^* \simeq C_4$. 

Suppose $u \in N(z)$ for some $u \in N(v)$ and $z \in V(H^*)$ with $d_{H^*}(z) = 2$. Since $\Delta(G) = 3$, $N[z] = \{u\} \cup (V(H^*) \setminus \{z'\})$ for some $z' \in V(H^*)$. If $G - N[z]$ is connected, then since $|V(G - N[z])| = 11$, the result follows as in Case~1 or by Lemmas~\ref{LemIsol1} and~\ref{ThmG8-11}. If $G - N[z]$ is not connected, then since $H_1$ is linked to $x$ and $x'$, we have ${\rm C}(G - N[z]) = \{G - (V(H^*) \cup \{u\}), (\{z'\}, \emptyset)\}$, so $\iota(G, P_3) \leq 1 + \iota(G - N[z], P_3)\leq 3$ by Lemmas~\ref{LemIsol1}, \ref{LemIsol2} and~\ref{ThmG8-11}. 

Now suppose $N(z) \cap N(v) = \emptyset$ for each $z \in V(H^*)$ with $d_{H^*}(z) = 2$. Then, since $H^*$ is linked to some $u_1 \in N(v)$, $u_1z_1 \in E(G)$ for some $z_1 \in V(H^*)$ with $d_{H^*}(z_1) = 1$. Thus, $H^* \simeq P_4$, meaning that $E(H^*) = \{z_1z_2, z_2z_3, z_3z_4\}$ for some $z_2, z_3, z_4 \in V(H^*) \setminus \{z_1\}$. For each $i \in \{2, 3\}$, since $d_{H^*}(z_i) = 2$, we have $N(z_i) \cap N(v) = \emptyset$, so $N(z_i) = \{z_{i-1}, z_{i+1}\}$. Let $G^* = G - (V(H^*) \cup \{u_1\})$ and $J = H^* - N_{H^*}[z_1]$. Since $|V(J)| = 2$, $\iota(J, P_3) = 0$. Since $H_1$ is linked to $x$ and $x'$, $G^*$ is connected. Since $|V(G^*)| = 10$, $\iota(G^*, P_3) \leq 2$ by Lemma~\ref{ThmG8-11}.

Suppose that $H^*$ is only linked to $u_1$. Then, ${\rm C}(G - N[z_1]) = \{G^*, J\}$. By Lemmas~\ref{LemIsol1} and~\ref{LemIsol2}, $\iota(G, P_3) \leq 1 + \iota(G - N[z_1], P_3) \leq 3$.

Now suppose that $H^*$ is linked to some $u_2 \in N(v) \setminus \{u_1\}$. Since $d(v) = 3$ and $\{u_1, u_2\}, \{x, x'\} \subseteq N(v)$, $\{u_1, u_2\} \cap \{x, x'\} \neq \emptyset$. We may assume that $u_1 = x$. Thus, $N(x) = \{v, y_1, z_1\}$. If $G - N[x]$ is connected, then the result follows as in Case~1 or by Lemmas~\ref{LemIsol1} and~\ref{ThmG8-11}. Suppose that $G - N[x]$ is not connected. 

Suppose $N(z_4) \cap N(v) = \emptyset$. Then, $H^* - z_1 \in {\rm C}(G - z_1)$. Let $X = \{x, z_1, z_2\}$. We have ${\rm C}(G - X) = \{G^*, J\}$, so $\iota(G - X, P_3) \leq 2$ by Lemma~\ref{LemIsol2}. Since $X \subseteq N[z_1]$, Lemma~\ref{LemIsol1} gives us $\iota(G, P_3) \leq 1 + \iota(G - X, P_3) \leq 3$. 

Now suppose $N(z_4) \cap N(v) \neq \emptyset$. Then, since $N(x) = \{v, y_1, z_1\}$, $N(z_4) \cap \{x', w\} \neq \emptyset$. Since $G - N[x]$ is not connected, either $N(z_4) \cap \{x', w\} = \{x'\}$, ${\rm C}(G - N[x]) = \{G - (N[x] \cup \{w\}), (\{w\}, \emptyset)\}$, and hence $\iota(G, P_3) \leq 1 + \iota(G - N[x], P_3) \leq 3$ by Lemmas~\ref{LemIsol1}, \ref{LemIsol2} and~\ref{ThmG8-11}, or $N(z_4) \cap \{x', w\} = \{w\}$ and $N(w) \setminus \{v, z_4\} \subseteq \{z_1\}$. Suppose that the latter holds. Then, $N(z_4) = \{w, z_3\}$. If $z_1 \in N(w)$, then ${\rm C}(G - N[w]) = \{G[V(H_1) \cup \{x, x'\}], G[\{z_2, z_3\}]\}$, so $\iota(G, P_3) \leq 1 + \iota(G - N[w], P_3) \leq 3$ by Lemmas~\ref{LemIsol1}, \ref{LemIsol2} and~\ref{ThmG8-11}. Suppose $z_1 \notin N(w)$. Then, $N(w) = \{v, z_4\}$ and ${\rm C}(G - N[x]) = \{A, B\}$ with $A = G[\{x'\} \cup (V(H_1) \setminus \{y_1\})]$ and $B = G[\{z_2, z_3, z_4, w\}]$. 
Let $L = \{y_2y_3, y_3y_4, \dots, y_6y_7\} \cup \{x'y'\}$. We have $L \subseteq E(A) \subseteq L \cup \{x'y''\}$ for some $y'' \in \{y_2, y_3, \dots, y_7\}$. Suppose $A \not\simeq C_7$. Since $|E(A)| \leq 7$, $A$ is not a $\mathcal{G}_7$-graph. By Lemma~\ref{ThmG7}, $\{v_A\} \in {\rm Is}(A, P_3)$ for some $v_A \in V(A)$, so $\{x, v_A, z_4\} \in {\rm Is}(G, P_3)$. Now suppose $A \simeq C_7$. Then, $N(x') = \{v, y', y''\} = \{v, y_2, y_7\}$, so $\{x', y_5, z_2\} \in {\rm Is}(G, P_3)$.\medskip

\noindent\textbf{Case~3:} \emph{$H_1$ is a $\mathcal{G}_3$-graph.} If $\sum_{H \in \mathcal{H} \setminus \mathcal{H}'} \iota(H,P_3) \leq 1$, then by Lemmas~\ref{LemIsol1} and~\ref{LemIsol2}, 
$$\iota(G, P_3) \leq 1 + \iota(G', P_3) = 1 + \iota(H_1, P_3) + \sum_{H\in\mathcal{H} \setminus \mathcal{H}'}\iota(H,P_3) \leq 3.$$ 
Suppose $\sum_{H \in \mathcal{H} \setminus \mathcal{H}'} \iota(H,P_3) \geq 2$. Recall that no member of $\mathcal{H} \setminus \mathcal{H}'$ is of order $3$. Since $\sum_{H \in \mathcal{H} \setminus \mathcal{H}'}|V(H)| = 8$ and no member of $\mathcal{H} \setminus \mathcal{H}'$ is an $\mathcal{E}$-graph, it follows by Lemma~\ref{ThmG7} that $\mathcal{H} \setminus \mathcal{H}'$ consists of either one $8$-vertex graph $I_1$ or two $4$-vertex graphs $I_1$ and $I_2$. Let $h = |\mathcal{H} \setminus \mathcal{H}'|$. Thus, $h \in [2]$. 

Suppose that for some $j \in [h]$ and $u \in N(v)$, $I_j$ is linked only to $u$. We may assume that $j = 1$. Let $G_1^* = G[V(H_1) \cup (N[v] \setminus \{u\})]$ and $G_2^* = G - (\{u\} \cup V(I_1))$. Then, $|V(G_1^*)| = 6$, and $G_1^*$ is connected as $H_1$ is linked to $x$ and $x'$. By Lemma~\ref{ThmG7}, since $G_1^* \not\simeq C_6$, $G_1^*$ has a $P_3$-isolating set $D_1$ with $|D_1| = 1$. Suppose $h = 2$. Then, $|V(G_2^*)| = 10$. If $I_2$ is also linked only to $u$, then $\{y_{u,I_1}, y_{u,I_2}\} \cup D_1$ is a $P_3$-isolating set of $G$ of size $3$. Suppose that $I_2$ is linked to some $u' \in N(v) \setminus \{u\}$. Then, $G_2^*$ is connected. By Lemma~\ref{ThmG8-11}, $G_2^*$ has a $P_3$-isolating set $D_2$ with $|D_2| \leq 2$. Thus, $\{y_{u,I_1}\} \cup D_2$ is a $P_3$-isolating set of $G$ of size at most $3$. Now suppose $h = 1$. Since $I_1$ is connected, $z \in N(y_{u,I_1})$ for some $z \in V(I_1)$. Let $I_1' = I_1 - N_{I_1}[y_{u,I_1}]$. Then, $|V(I_1')| \leq 6$. Let $D_3$ be a smallest $P_3$-isolating set of $I_1'$. If $|D_3| \leq 1$, then $\{y_{u,I_1}\} \cup D_1 \cup D_3$ is a $P_3$-isolating set of $G$ of size at most $3$. Suppose $|D_3| \geq 2$. By Lemma~\ref{ThmG7}, we clearly obtain ${\rm C}(I_1') = \{T_1, T_2\}$ for some distinct $\mathcal{G}_3$-graphs $T_1$ and $T_2$. Since $|V(I_1)| = 8$, $z$ is the only neighbour of $y_{u, I_1}$ in $I_1$. Since $I_1$ is connected, $z_1, z_2 \in N(z)$ for some $z_1 \in V(T_1)$ and $z_2 \in V(T_2)$. This yields $\{z\} \in {\rm Is}(I_1, P_3)$, contradicting $\sum_{H \in \mathcal{H} \setminus \mathcal{H}'} \iota(H,P_3) \geq 2$.

Now suppose that for each $j \in [h]$, there exist $u_j^1, u_j^2 \in N(v)$ such that $u_j^1 \neq u_j^2$ and $I_j$ is linked to $u_j^1$ and $u_j^2$. If $u \in N(z)$ for some $u \in N(v)$ and $z \in V(H_1)$ with $d_{H_1}(z) = 2$, then $N[z] = \{u\} \cup V(H_1)$ (as $\Delta(G) = 3$), $G - N[z]$ is connected (as $I_j$ is linked to $u_j^1$ and $u_j^2$), $|V(G - N[z])| = 11$, and hence the result follows as in Case~1 or by Lemmas~\ref{LemIsol1} and~\ref{ThmG8-11}. Suppose $N(z) \cap N(v) = \emptyset$ for each $z \in V(H_1)$ with $d_{H_1}(z) = 2$. Then, since $H_1$ is linked to $x$, $xy_1 \in E(G)$ for some $y_1 \in V(H_1)$ with $d_{H_1}(y_1) = 1$. Thus, $H_1 \simeq P_3$, and hence $E(H_1) = \{y_1y_2, y_2y_3\}$ for some $y_2, y_3 \in V(H_1) \setminus \{y_1\}$, and $N(y_2) = \{y_1, y_3\}$. 

Suppose $N(y_3) \cap N(v) \subseteq \{x\}$. Since $H_1$ is linked to $x'$, $N(y_1) = \{x, x', y_2\}$. Let $Y = \bigcup_{i=1}^h V(I_i)$. Let $G_Y = G[\{v, w\} \cup Y]$. Thus, $|V(G_Y)| = 10$, and $\{y_1\} \cup D_Y \in {\rm Is}(G, P_3)$ for any $D_Y  \in {\rm Is}(G_Y, P_3)$. If $G_Y$ is connected, then $\iota(G_Y, P_3) \leq 2$ by Lemma~\ref{ThmG8-11}, so $\iota(G, P_3) \leq 3$. Suppose that $G_Y$ is not connected. Then, for some $j \in [h]$, $I_j$ is not linked to $w$. We may assume that $j = 1$. If $h = 1$, then ${\rm C}(G_Y) = \{G[\{v, w\}], I_1\}$, $\iota(G_Y, P_3) = \iota(I_1, P_3) \leq 2$ by Lemmas~\ref{LemIsol2} and~\ref{ThmG8-11}, and hence $\iota(G, P_3) \leq 3$. Suppose $h = 2$. Then, ${\rm C}(G_Y)$ is $\{G[\{v, w\}], I_1, I_2\}$ or $\{G[\{v, w\} \cup V(I_2)], I_1\}$. By Lemmas~\ref{LemIsol2} and~\ref{ThmG7}, $\iota(G_Y, P_3) \leq 2$, so $\iota(G, P_3) \leq 3$. 

Now suppose $N(y_3) \cap N(v) \nsubseteq \{x\}$. We may assume that $x' \in N(y_3)$. Since $d(v) = 3$ and $\{u_1^1, u_1^2\}, \{x, x'\} \subseteq N(v)$, $\{u_1^1, u_1^2\} \cap \{x, x'\} \neq \emptyset$. We may assume that $u_1^1 = x$. Thus, $N(x) = \{v, y_1, y_{x,I_1}\}$. Since $G[\{v, x, y_1, y_2, y_3, x'\}] \not\simeq C_6$, we obtain $x'y_1 \in E(G)$, so $N(x') = \{v, y_1, y_3\}$ and $u_1^2 = w$. Suppose $h = 2$. Then, we have $u_2^1, u_2^2 \in N(v) \setminus \{x, x'\} = \{w\}$, contradicting $u_2^1 \neq u_2^2$. Thus, $h = 1$. Let $G^* = G - (\{v, x'\} \cup V(H_1))$. Since $xy_{x,I_1}, wy_{w,I_1} \in E(G)$, $G^*$ is connected. By Lemma~\ref{ThmG8-11}, since $|V(G^*)| = 10$, $G^*$ has a $P_3$-isolating set $D^*$ with $|D^*| \leq 2$. Since $\{x'\} \cup D^* \in {\rm Is}(G, P_3)$, $\iota(G, P_3) \leq 3$.\qed
}
\\

\noindent
\textbf{Proof of Lemma~\ref{ThmP3small}.} The result is an immediate consequence of Lemmas~\ref{ThmG7}--\ref{ThmG15}. 
\qed

\section{Acknowledgements}
The authors are grateful to the anonymous referees for checking the paper and providing constructive remarks. They would also like to thank Juan Scerri for computational assistance in confirming the graphs in $\mathcal{G}_{11}$. The research work disclosed in this publication is partially funded by the Tertiary Education Scholarships Scheme (Malta).


\begin{thebibliography}{}

\bibitem{BBS} K. Bartolo, P. Borg, D. Scicluna, Isolation of squares in graphs, Discrete Math. 347 (2024), paper 114161.

\bibitem{Borg1} P. Borg, Isolation of cycles, Graphs Combin. 36 (2020), 631--637. 

\bibitem{BorgIsolConnected2023} P. Borg, Isolation of connected graphs, Discrete Appl. Math. 339 (2023), 154--165.

\bibitem{Borgisdom} P. Borg, Proof of a conjecture on isolation of graphs dominated by a vertex, Discrete Appl. Math. 371 (2025), 247-253.

\bibitem{Borgrsc} P. Borg, Isolation of regular graphs, stars and $k$-chromatic graphs, Discrete Math. 349 (2026), paper 114706.

\bibitem{BK} P. Borg, P. Kaemawichanurat, Partial domination of maximal outerplanar graphs, Discrete Appl. Math. 283 (2020), 306--314.

\bibitem{BK2} P. Borg, P. Kaemawichanurat, Extensions of the Art Gallery Theorem, Ann. Comb. 27 (2023), 31--50.

\bibitem{BFK} P. Borg, K. Fenech, P. Kaemawichanurat, Isolation of $k$-cliques, Discrete Math. 343 (2020), paper 111879.

\bibitem{BFK2} P. Borg, K. Fenech, P. Kaemawichanurat, Isolation of $k$-cliques II, Discrete Math. 345 (2022), paper 112641.

\bibitem{BG} G. Boyer, W. Goddard, Disjoint isolating sets and graphs with maximum isolation number, Discrete Appl. Math. 356 (2024), 110--116.

\bibitem{CaWa13} C.N. Campos, Y. Wakabayashi, On dominating sets of maximal outerplanar graphs, Discrete Appl. Math. 161 (2013), 330--335.

\bibitem{CaHa17} Y. Caro, A. Hansberg, Partial domination - the isolation number of a graph, Filomat 31 (2017), 3925--3944.  

\bibitem{CLWX} J. Chen, Y. Liang, C. Wang, S. Xu, Algorithmic aspects of $\{P_k\}$-isolation in graphs and extremal graphs for a $\{P_3\}$-isolation bound, Inf. Process. Lett. 187 (2025), paper 106521.

\bibitem{Ch75} V. Chv\'{a}tal, A combinatorial theorem in plane geometry, J. Combin. Theory Ser. B 18 (1975), 39--41.

\bibitem{C} E.J. Cockayne, Domination of undirected graphs -- A survey, Lecture Notes in Mathematics, vol. 642, Springer, 1978, 141--147.

\bibitem{CH} E.J. Cockayne, S.T. Hedetniemi, Towards a theory of domination in graphs, Networks 7 (1977), 247--261.
 
\bibitem{CZZ} Q. Cui, J. Zhang, L. Zhong, Extremal graphs for the $K_{1,2}$-isolation number of graphs, Bull. Malays. Math. Sci. Soc. 47 (2024), paper 115.

\bibitem{DoHaJo16} M. Dorfling, J.H. Hattingh, E. Jonck, Total domination in maximal outerplanar graphs II, Discrete Appl. Math. 339 (2016), 1180--1188.

\bibitem{DoHaJo17} M. Dorfling, J.H. Hattingh, E. Jonck, Total domination in maximal outerplanar graphs, Discrete Appl. Math. 217 (2017), 506--511.

\bibitem{HHS} T.W. Haynes, S.T. Hedetniemi, P.J. Slater, Fundamentals of Domination in Graphs, Marcel Dekker, Inc., New York, 1998.

\bibitem{HHS2} T.W. Haynes, S.T. Hedetniemi, P.J. Slater (Editors), Domination in Graphs: Advanced Topics, Marcel Dekker, Inc., New York, 1998.

\bibitem{HL2} S.T. Hedetniemi, R.C. Laskar, Bibliography on domination in graphs and some basic definitions of domination parameters, Discrete Math. 86 (1990), 257--277.

\bibitem{HL} S.T. Hedetniemi, R.C. Laskar (Editors), Topics on Domination, in: Annals of Discrete Mathematics, vol. 48, North-Holland Publishing Co., Amsterdam, 1991, Reprint of Discrete Math. 86 (1990), no. 1--3.

\bibitem{HeKa18} M.A. Henning, P. Kaemawichanurat, Semipaired domination in maximal outerplanar graphs, J. Comb. Optim. 38 (2019), 911--926.

\bibitem{HZJ} Y. Huang, G. Zhang, X. Jin, New results on the $1$-isolation number of graphs without short cycles, Discrete Appl. Math. 379 (2026), 222--235. 

\bibitem{LeZuZy17} M. Lema\'{n}ska, R. Zuazua, P. \.{Z}yli\'{n}ski, Total dominating sets in maximal outerplanar graphs, Graphs Combin. 33 (2017), 991--998.

\bibitem{LMS} M. Lema\'{n}ska, M. Mora, M.J. Souto--Salorio, Graphs with isolation number equal to one third of the order, Discrete Math. 347 (2024), paper 113903.

\bibitem{Li16} Z. Li, E. Zhu, Z. Shao, J. Xu, On dominating sets of maximal outerplanar and planar graphs, Discrete Appl. Math. 198 (2016), 164--169.

\bibitem{MaTa96} L.R. Matheson, R.E. Tarjan, Dominating sets in planar graphs, European J. of Combin. 17 (1996), 565--568.

\bibitem{McKay} B. McKay, \url{https://users.cecs.anu.edu.au/~bdm/data/graphs.html}, accessed 30/12/2024.

\bibitem{Ore} O. Ore, Theory of graphs, American Mathematical Society Colloquium Publications, vol. 38, American Mathematical Society, Providence, R.I., 1962. 

\bibitem{To13} S. Tokunaga, Dominating sets of maximal outerplanar graphs, Discrete Appl. Math. 161 (2013), 3097--3099.

\bibitem{KaJi} S. Tokunaga, T. Jiarasuksakun, P. Kaemawichanurat, Isolation number of maximal outerplanar graphs, Discrete Appl. Math. 267 (2019), 215--218.

\bibitem{West} D.B. West, Introduction to graph theory, second ed., Prentice Hall, 1996.

\bibitem{ZW} G. Zhang, B. Wu, $K_{1,2}$-isolation in graphs, Discrete Appl. Math. 304 (2021), 365--374. 

\bibitem{Z} P. \.{Z}yli\'{n}ski, Vertex-edge domination in graphs, Aequat. Math. 93 (2019), 735--742.

\end{thebibliography}
\end{document}